\documentclass[11pt]{amsart}
\usepackage{amsthm, amsmath}
\usepackage{amssymb} 
\usepackage{graphicx}
\usepackage{epsf}
\usepackage[all]{xy}

\theoremstyle{plain}
\newtheorem{theorem}{Theorem}[section]
\newtheorem{lemma}[theorem]{Lemma}
\newtheorem{cor}[theorem]{Corollary}
\newtheorem{prop}[theorem]{Proposition}
\newtheorem{conj}[theorem]{Conjecture}
\newtheorem{definition}[theorem]{Definition}

\numberwithin{equation}{subsection}

\theoremstyle{remark}

\theoremstyle{definition}

\newtheorem{Remark}[theorem]{Remark}
\newtheorem{Question}[theorem]{Question}

\newcommand{\G}{{\mathbb G}_m}

\newcommand\xt{[\hspace{-.12em}[ t ]\hspace{-.12em} ] } 
\newcommand\xT{(\!( t )\!)}
\newcommand\cQ{{\mathcal Q}}

\title{Equidimensionality of convolution morphisms and applications to saturation problems}
\author{Thomas J. Haines}
\thanks{Research supported by NSF grant DMS 0303605, and by a Sloan Research Fellowship.}

\begin{document}

\begin{abstract}
Fix a split connected reductive group $G$ over a field $k$, and a positive integer $r$.  For any $r$-tuple of dominant coweights $\mu_i$ of $G$, we consider the restriction $m_{\mu_\bullet}$ of the r-fold convolution morphism of Mirkovic-Vilonen \cite{MV1,MV2} to the twisted product of affine Schubert varieties corresponding to $\mu_\bullet$.  We show that if all the coweights $\mu_i$ are minuscule, then the fibers of $m_{\mu_\bullet}$ are equidimensional varieties, with dimension the largest allowed by the semi-smallness of $m_{\mu_\bullet}$.  We derive various consequences: the equivalence of the non-vanishing of Hecke and representation ring structure constants, and a saturation property for these structure constants, when the coweights $\mu_i$ are sums of minuscule coweights.  This complements the saturation results of Knutson-Tao \cite{KT} and Kapovich-Leeb-Millson \cite{KLM}.  We give a new proof of the P-R-V conjecture in the ``sums of minuscules'' setting.  Finally, we generalize and reprove a result of Spaltenstein pertaining to equidimensionality of certain partial Springer resolutions of the nilpotent cone for ${\rm GL}_n$.
\end{abstract}

\maketitle

\markboth{THOMAS J. HAINES}{EQUIDIMENSIONALITY OF CERTAIN CONVOLUTION MORPHISMS}

\bigskip
\section{Introduction}

Let $G$ be a split connected reductive group over a finite field $\mathbb F_q$, with Langlands dual $\widehat{G} = \widehat{G}(\overline{\mathbb Q}_\ell)$, where ${\rm char}(\mathbb F_q) = p$ and $\ell \neq p$ is prime.  The geometric Satake isomorphism of Mirkovic-Vilonen \cite{MV2} establishes a geometric construction of $\widehat{G}$.  More precisely, it identifies $\widehat{G}$ with the automorphism group of the fiber functor of a certain Tannakian category.  Letting $F = \mathbb F_q \xT$ and $\mathcal O = \mathbb F_q \xt$, the latter is the category ${\rm P}_{G(\mathcal O)}$ of $G(\mathcal O)$-equivariant perverse $\overline{\mathbb Q}_{\ell}$-sheaves ${\mathcal F}$ on the affine Grassmannian
$$
{\mathcal Q} = G(F)/G(\mathcal O),
$$
viewed as an ind-scheme over $\mathbb F_q$.  The fiber functor 
$$
\mathcal F \mapsto {\rm H}^*({\mathcal Q}, \mathcal F)
$$
takes ${\rm P}_{G(\mathcal O)}$ to the category of graded finite-dimensional $\overline{\mathbb Q}_\ell$-vector spaces.  In order to give ${\rm P}_{G(\mathcal O)}$ a Tannakian structure, one needs to endow it with a tensor product with commutativity and associativity constraints.  There are a few different ways to construct the tensor product (see especially \cite{Gi}, \cite{MV1}, and \cite{Ga}).  The present article will use the construction in \cite{MV1}, which is defined in terms of the {\em convolution morphism} 
$$
m_{\mu_\bullet}: \overline{\mathcal Q}_{\mu_1} \widetilde{\times} \cdots \widetilde{\times} \overline{\mathcal Q}_{\mu_r} \rightarrow \overline{\mathcal Q}_{|\mu_\bullet|}.
$$
Here the $\mu_i$ are dominant cocharacters of $G$ indexing various $G(\mathcal O)$-orbits ${\mathcal Q}_{\mu_i} \subset {\mathcal Q}$ (via the Cartan decomposition), $|\mu_\bullet| := \sum_i \mu_i$, and the morphism 
$m_{\mu_\bullet}$ forgets all but the last element in the twisted product (see section \ref{notation}).  The morphism 
$m_{\mu_\bullet}$ is used to construct the $r$-fold convolution product in ${\rm P}_{G(\mathcal O)}$, as follows.  
Given $G(\mathcal O)$-equivariant perverse sheaves ${\mathcal F}_1, \dots, {\mathcal F}_r$, supported on various closures $\overline{\mathcal Q}_{\mu_1}, \dots, \overline{\mathcal Q}_{\mu_r}$, there is a well-defined perverse ``twisted external product'' sheaf ${\mathcal F}_1 \widetilde{\boxtimes} \cdots \widetilde{\boxtimes} {\mathcal F}_{r}$ on the twisted product $\overline{\mathcal Q}_{\mu_1} \widetilde{\times} \cdots \widetilde{\times} \overline{\mathcal Q}_{\mu_r}$; see section \ref{notation}.  Then the $r$-fold convolution product is defined by the proper push-forward on derived categories
$$
{\mathcal F}_1 * \cdots * {\mathcal F}_r = m_!({\mathcal F}_1 \widetilde{\boxtimes} \cdots \widetilde{\boxtimes} {\mathcal F}_{r}).
$$  

For brevity, let us write $K = G(\mathcal O)$, a maximal compact subgroup of the loop group $G(F)$. 
Zariski-locally the twisted product $\overline{\mathcal Q}_{\mu_1} \widetilde{\times} \cdots \widetilde{\times} \overline{\mathcal Q}_{\mu_r}$ is just the usual product and the morphism $m_{\mu_\bullet}$ is given by
$$
m_{\mu_\bullet}: (g_1 K, g_2 K , \dots, g_r K) \mapsto g_1g_2 \cdots g_r \, K.
$$
Using this one may check that 
under the sheaf-function dictionary \`{a} la Grothendieck, the tensor structure on ${\rm P}_{G(\mathcal O)}$ corresponds to the usual convolution in the spherical Hecke algebra ${\mathcal H}_q = C_c(K \backslash G(F)/K)$.  This is the convolution algebra of compactly-supported $\overline{\mathbb Q}_\ell$-valued functions on $G(F)$ which are bi-invariant under $K$, where the convolution product (also denoted $*$) is defined using the Haar measure which gives $K$ volume 1.   This is the reason why we call $m_{\mu_\bullet}$ a {\em convolution morphism}.  

The morphism $m_{\mu_\bullet}$ is projective, birational, and semi-small and locally-trivial in the stratified sense; see \cite{MV1}, \cite{NP} and $\S \ref{local_triviality}$ for proofs of these properties, and \cite{H} for some further discussion.  These properties  are essential for the construction of the tensor product on ${\rm P}_{G(\mathcal O)}$.

As is well-known, the fibers of the morphism $m_{\mu_\bullet}$ carry representation-theoretic information (see section \ref{content_of_fibers}).  The purpose of this article is to establish a new equidimensionality property of these fibers in a very special situation, and then to extract some consequences of combinatorial and representation-theoretic nature.  The main result is the following theorem.  Let $\rho$ denote the half-sum of the positive roots for $G$, and recall that the semi-smallness of $m_{\mu_\bullet}$  means that for every $y \in {\mathcal Q}_\lambda \subset \overline{\mathcal Q}_{|\mu_\bullet|}$, the fiber over $y$ satisfies the following bound on its dimension
$$
{\rm dim}(m_{\mu_\bullet}^{-1}(y)) \leq \frac{1}{2}[{\rm dim}(\overline{\mathcal Q}_{|\mu_\bullet|}) - {\rm dim}(\overline{\mathcal Q}_\lambda)] = \langle \rho, |\mu_\bullet| - \lambda \rangle.
$$

\begin{theorem} [Equidimensionality for minuscule convolutions]  \label{thm_A} Let $y \in {\mathcal Q}_\lambda \subset \overline{\mathcal Q}_{|\mu_\bullet|}$.  Suppose each coweight $\mu_i$ is minuscule.  
Then every irreducible component of the fiber $m_{\mu_\bullet}^{-1}(y)$ has dimension $\langle \rho, |\mu_\bullet| - \lambda \rangle$.  
\end{theorem}
Recall that a coweight $\mu$ is {\em minuscule} if $\langle \alpha, \mu \rangle \in \{ -1, 0, 1\}$ for every root $\alpha$.  The following result is a corollary of the proof.

\begin{cor} \label{paved_by_affines} If every $\mu_i$ is minuscule, then each fiber $m_{\mu_\bullet}^{-1}(y)$ admits a paving by affine spaces.
\end{cor}

The conclusions in Theorem \ref{thm_A} fail without the hypothesis that each $\mu_i$ is minuscule. 
Without that hypothesis, the dimension of the fiber can be strictly less than $\langle \rho, |\mu_\bullet| - \lambda \rangle$.  This can happen even if we weaken the hypothesis to ``each $\mu_i$ is minuscule or quasi-minuscule'', see Remark \ref{quasi_min_remark}.  Further, even for $G = {\rm GL}_n$ there exist coweights of the form $\mu_i = (d_i,0^{n-1})$ where $d_1 + \cdots + d_r = n$, for which certain fibers $m_{\mu_\bullet}^{-1}(y)$ are {\em not} equidimensional, see Remark \ref{nonequidim}.  We do not know how to characterize the tuples $\mu_\bullet$ for which every fiber $m_{\mu_\bullet}^{-1}(y)$ is paved by affine spaces, see Question \ref{questions}.

\medskip

Nevertheless, a similar equidimensionality statement continues to hold when we require each $\mu_i$ to be a sum of minuscules (see $\S \ref{further_equidim_section}$).  In its most useful form it concerns the intersection of the fiber $m_{\mu_\bullet}^{-1}(y)$ with the open stratum ${\mathcal Q}_{\mu_\bullet} = {\mathcal Q}_{\mu_1} 
\widetilde{\times} \cdots \widetilde{\times} {\mathcal Q}_{\mu_r}$ of the twisted product $\widetilde{\mathcal Q}_{\mu_\bullet}$.  The following result is an easy corollary of Theorem \ref{thm_A}.  It is proved in Proposition \ref{further_equidim} (see also \cite{H}, $\S 8$).

\begin{theorem} [Equidimensionality for sums of minuscules] \label{thm_B}
Suppose each $\mu_i$ is a sum of dominant minuscule coweights.  Then the intersection
$$
m_{\mu_\bullet}^{-1}(y) \cap {\mathcal Q}_{\mu_\bullet}
$$
is equidimensional of dimension $\langle \rho, |\mu_\bullet| - \lambda \rangle$, provided the intersection is non-empty.
\end{theorem}

This result also generally fails to hold without the hypothesis on the coweights $\mu_i$ (see Remark 
\ref{quasi_min_remark}).  Note that Theorem \ref{thm_A} is actually a special case of 
Theorem \ref{thm_B}.

\smallskip

Theorem \ref{thm_B} allows us to establish a relation between structure constants of Hecke and representation rings, generalizing \cite{H}, which treated the case of ${\rm GL}_n$.  Namely, thinking of $(\mu_\bullet,\lambda)$ as an $r+1$-tuple of dominant weights of $\widehat{G}$ (resp. coweights of $G$), we may define structure constants ${\rm dim}(V^\lambda_{\mu_\bullet})$  (resp. $c^\lambda_{\mu_\bullet}(q)$) for the representation 
ring of the category ${\rm Rep}(\widehat{G})$ (resp. for the Hecke algebra ${\mathcal H}_q$) corresponding to the multiplication of basis elements consisting of highest-weight representations $V_{\mu_1}, \dots, V_{\mu_r}$ (resp. characteristic functions $f_{\mu_1} = 1_{K\mu_1K}, \dots, f_{\mu_r} = 1_{K\mu_rK}$).  In other words, we consider the decompositions
\begin{align*}
V_{\mu_1} \otimes \cdots \otimes V_{\mu_r} ~ &= ~ \bigoplus_\lambda  \, V^\lambda_{\mu_\bullet} \otimes V_\lambda \\
f_{\mu_1} * \cdots * f_{\mu_r} ~ &= ~ \sum_\lambda c^\lambda_{\mu_\bullet}(q) \, f_\lambda
\end{align*}
in ${\rm Rep}(\widehat{G})$ and ${\mathcal H}_q$, respectively. 
Following \cite{H}, consider the properties
\begin{align*}
{\rm Rep}(\mu_\bullet,\lambda) ~ &: ~ {\rm dim}(V^\lambda_{\mu_\bullet}) > 0 \\
{\rm Hecke}(\mu_\bullet, \lambda) ~ &: ~ c^\lambda_{\mu_\bullet}(q) \neq 0.
\end{align*}

It is a general fact that ${\rm Rep}(\mu_\bullet, \lambda) \Rightarrow {\rm Hecke}(\mu_\bullet, \lambda)$, for all groups $G$ (see \cite{KLM}, Theorem 1.13, and Corollary \ref{structure_constants_formula} below).  The reverse implication holds for ${\rm GL}_n$, but fails for general tuples $\mu_\bullet$ attached to other groups (see \cite{KLM},\cite{H}, and Remark \ref{quasi_min_remark}).  The following consequence of Theorem \ref{thm_B} shows that there is a 
natural condition on the coweights $\mu_i$ which ensures that the reverse implication does hold. 

\begin{theorem} [Equivalence of non-vanishing of structure constants] \label{thm_C}
If each $\mu_i$ is a sum of dominant minuscule coweights of $G$, then
$$
{\rm Rep}^{\widehat{G}}(\mu_\bullet, \lambda) \Leftrightarrow {\rm Hecke}^G(\mu_\bullet, \lambda).
$$
\end{theorem}

 Since every coweight of ${\rm GL}_n$ is a sum of minuscule coweights, this puts the ${\rm GL}_n$ case into a broader context. 
For groups not of type $A$, many (or all) coweights are not sums of minuscules, and this is reflected by the abundance of counterexamples to the implication ${\rm Hecke}(\mu_\bullet, \lambda) \Rightarrow {\rm Rep}(\mu_\bullet,\lambda)$ for those groups.

\smallskip

As first pointed out by M. Kapovich, B. Leeb, and J. Millson \cite{KLM}, the translation from the representation ring structure constants to Hecke algebra structure constants has some applications, in particular to saturation questions for general groups.  The authors of \cite{KLM} investigated saturation questions for the structure constants of $\mathcal H_q$, and their results apply to general groups $G$.  Results such as Theorem \ref{thm_C} allow us to deduce saturation theorems for ${\rm Rep}(\widehat{G})$.  

\begin{theorem} [A saturation theorem for sums of minuscules] \label{thm_D}
Suppose $\mu_\bullet$ is an $r$-tuple of dominant weights for $\widehat{G}$, whose sum belongs to the root lattice of $\widehat{G}$.  Suppose each $\mu_i$ is a sum of dominant minuscule weights.  Let $V_{\mu_i}$ denote the irreducible $\widehat{G}$-module with highest weight $\mu_i$.  Then
\begin{enumerate}
\item[(1)] If $k= k_G$ denotes the Hecke algebra saturation factor for $G$ as defined in \cite{KLM}, then 
$$
(V_{N\mu_1} \otimes \cdots \otimes V_{N\mu_r})^{\widehat{G}} \neq 0 \Rightarrow (V_{k\mu_1} \otimes \cdots \otimes V_{k\mu_r})^{\widehat{G}} \neq 0,
$$
for every positive integer $N$.
\item[(2)]  If the simple factors of $G_{ad}$ are all of type $A, B, C$ or $E_7$, then the above implication holds with $k$ replaced by 1.
\end{enumerate}
\end{theorem}

The analogue of part (1) for Hecke algebra structure constants is due to M. Kapovich, B. Leeb, and J. Millson \cite{KLM}.  We derive part (1) from their result by applying Theorem \ref{thm_C}, with $\lambda = 0$.   In fact a sharper version of part (1) is valid: we need only assume that at least $r-1$ of the weights $\mu_i$ are sums of minuscules, see Theorem \ref{saturation_for_sums_of_minuscules}.

A somewhat more comprehensive version of part (2) is proved in Theorem \ref{evidence_for_conjecture}, again by establishing the Hecke algebra analogue.  That analogue is proved in Theorem \ref{saturation_for_Hecke_constants} of the Appendix, written jointly with M. Kapovich and J. Millson.  Based on this result and some computer calculations done using LiE, we conjecture that the conclusion of part (2) holds in all cases (i.e. factors of type $D$ and $E_6$ should also be allowed; see Conjecture \ref{conjecture}).  

Note that for $\widehat{G} = {\rm GL}_n(\mathbb C)$, part (2) is not new.  It is the well-known saturation property of ${\rm GL}_n$, which was first proved by A. Knutson and T. Tao in their paper \cite{KT}.  The Hecke algebra approach was introduced in \cite{KLM}, which provided a new proof of the Knutson-Tao result, and suggested that saturation problems for more general groups are best approached via Hecke algebras and triangles in Bruhat-Tits buildings. 

\smallskip

In their recent preprint \cite{KM}, Kapovich and Millson have announced some results which are closely related to our Theorems \ref{thm_C} and \ref{thm_D}, and which are proved by completely different methods; see Remarks \ref{no_assumption_on_lambda}, \ref{KM_refinement}.

\smallskip

Theorems \ref{thm_A}, \ref{thm_B}, and \ref{thm_C} were proved for ${\rm GL}_n$ in \cite{H}, as consequences of the geometric Satake isomorphism, the P-R-V property, and Spaltenstein's theorem in \cite{Sp} on the equidimensionality of certain partial Springer resolutions.  In this paper, the geometric Satake isomorphism (more precisely, a corollary of it, Theorem \ref{weak_GS}) remains a key ingredient, and in some sense this work could be viewed as an application of that powerful result.  On the other hand, the present proofs of Theorems \ref{thm_A}-\ref{thm_C} rely on neither the P-R-V property nor Spaltenstein's theorem.  In fact, here we turn the logic around, giving a new proof of the P-R-V property in the ``sums of minuscules'' situation, and also giving a new proof and a generalization of Spaltenstein's theorem.  Those results are explained in sections \ref{PRV_section}  and \ref{Springer_section}, respectively.

\medskip

\noindent {\em Acknowledgments.}  
I express my thanks to Misha Kapovich and John Millson for generously sharing their ideas with me, especially in relation to the Appendix, which was written jointly with them.  I also thank them for giving me early access to their recent work \cite{KM}.  I am indebted to Jeff Adams for his invaluable help with LiE; his programs were used to run some extensive computer checks of Conjecture \ref{conjecture}.  Finally, I thank the referee for some very helpful suggestions for simplifying the proof of Theorem \ref{equidim}.

\medskip

\section{Preliminaries and notation} \label{notation}

\subsection{General notation}

Let $k$ denote a field, usually taken to be the complex numbers $\mathbb C$, a finite field ${\mathbb F}_q$, or an algebraic closure $\overline{\mathbb F}_q$ of a finite field.  Let $\mathcal O = k \xt$ (resp. $F = k \xT$) denote the ring of formal power series (resp. Laurent series) over $k$.  

Let $G$ denote a split connected reductive group over $k$.  Fix a $k$-split maximal torus $T$ and a $k$-rational Borel subgroup $B$ containing $T$.  We have $B = TU$, where $U$ is the unipotent radical of $B$.  Let $X_+ \subset X_*(T)$ denote the set of $B$-dominant integral coweights for $G$.  By $W$ we denote the finite Weyl group $N_G(T)/T$.  The Bruhat order $\leq$ on $W$ will always be the one determined by the Borel $B$ we have fixed.  Let $w_0$ denote the longest element in $W$.  
 
Consider the ``loop group'' $G(F) = G(k\xT)$ as an ind-scheme
over $k$.  Occasionally we designate this by $LG$, and the ``maximal compact'' subgroup $G(\mathcal O)$ by $L^{\geq 0}G$ or simply $K$.   The affine Grassmannian $\cQ$ (over the field $k$) is the 
$fpqc$-quotient sheaf
 $G(k\xT)/G(k\xt)$; it is an ind-scheme. If $G=GL_n$ and $R$ is a $k$-algebra, $\cQ(R)$
is the set of all $R\xt$-lattices in $R\xT^n$. If $G=GSp_{2n}$, it is the set of
lattices in $R\xT^{2n}$ which are self-dual up to an element in $R\xt^\times$. 

By the Cartan decomposition we have a stratification into $G(k\xt)$-orbits:
$$ {\mathcal Q} = \coprod_{\mu \in X_+} G(k\xt)\mu G(k\xt)/G(k\xt). $$
Here we embed $X_*(T)$ into $G(k\xT)$ by the rule $\mu \mapsto \mu(t) \in T(k\xT)$.  We will denote the $G(k\xt)$-orbit of $\mu$ simply by $\cQ_\mu$ in the sequel.  The closure relations are determined by the standard partial order $\preceq$ on dominant coweights: $\cQ_\lambda \subset \overline{\cQ}_\mu$ if and only $\lambda \preceq \mu$, which by definition holds if and only if $\mu - \lambda$ is a sum of $B$-positive coroots. Given $L, L' \in {\mathcal Q}$, let ${\rm inv}(L,L') \in X_+$ denote the relative position of  $L,L'$, where by definition
$$
{\rm inv}(gK,g'K) = \lambda \Leftrightarrow g^{-1}g' \in K\lambda K.
$$

There is a canonical perfect pairing $\langle \cdot ~,~ \cdot \rangle : X^*(T) \times X_*(T) \rightarrow \mathbb Z$.   Let $\rho$ denote the half-sum of the $B$-positive roots of $G$.  Given $\mu \in X_+$, the $K$-orbit ${\mathcal Q}_\mu$ is a smooth quasiprojective variety of dimension $\langle 2\rho, \mu \rangle$ over $k$.  Let $\overline{\mathcal Q}_\mu \subset {\mathcal Q}$ denote the closure of ${\mathcal Q}_\mu$ in the ind-scheme ${\mathcal Q}$.  

Let $e_0$ denote the base point in the affine Grassmannian for $G$, i.e., the point corresponding to the coset $K \in G(F)/K$.  For $\nu \in X_*(T)$, let $t_\nu := \nu(t) \in LG$.  For a dominant coweight $\lambda$, denote $e_\lambda = t_\lambda e_0$.

Now let $\mu_\bullet = (\mu_1, \dots, \mu_r)$, where $\mu_i \in X_+$ for $1 \leq i \leq r$.  We define the twisted product scheme
$$
\widetilde{\mathcal Q}_{\mu_\bullet} = \overline{\mathcal Q}_{\mu_1} \widetilde{\times} \cdots \widetilde{\times} 
\overline{\mathcal Q}_{\mu_r}
$$
to be the subscheme of ${\mathcal Q}^r$ consisting of points $(L_1, \dots, L_r)$ such that ${\rm inv}(L_{i-1},L_i) \preceq \mu_i$ for $1 \leq i \leq r$ (letting $L_0 = e_0$).  The projection onto the last coordinate gives the proper surjective birational morphism
$$
m_{\mu_\bullet}: \widetilde{\mathcal Q}_{\mu_\bullet} \rightarrow \overline{\mathcal Q}_{|\mu_\bullet|},
$$
where by definition $|\mu_\bullet| = \sum_i \mu_i$.  

Note that the target of $m_{\mu_\bullet}$ is stratified by the $K$-orbits ${\mathcal Q}_\lambda$ for $\lambda$ ranging over dominant coweights satisfying $\lambda \preceq |\mu_\bullet|$.  Similarly, the domain is stratified by the locally closed twisted products ${\mathcal Q}_{\mu'_\bullet} := {\mathcal Q}_{\mu'_1} \widetilde{\times} \cdots \widetilde{\times} {\mathcal Q}_{\mu'_r}$, where $\mu'_i$ ranges over dominant coweights satisfying $\mu'_i \preceq \mu_i$.  Here ${\mathcal Q}_{\mu'_\bullet}$ is defined exactly as is $\widetilde{\mathcal Q}_{\mu'_\bullet}$, except that the conditions ${\rm inv}(L_{i-1},L_i) \preceq \mu'_i$ are replaced with ${\rm inv}(L_{i-1},L_i) = 
\mu'_i$.  

With respect to these stratifications, $m_{\mu_\bullet}$ is locally trivial and semi-small (in the stratified sense).  The local triviality is discussed in $\S \ref{local_triviality}$.  The semi-smallness means that for every inclusion ${\mathcal Q}_\lambda \subset m_{\mu_\bullet}({\mathcal Q}_{\mu'_\bullet})$, the fibers of the restricted morphism
$$
m_{\mu_\bullet}: m_{\mu_\bullet}^{-1}({\mathcal Q}_\lambda) \cap {\mathcal Q}_{\mu'_\bullet} \rightarrow {\mathcal Q}_\lambda
$$
have dimension bounded above by 
$$
\dfrac{1}{2}[{\rm dim}({\mathcal Q}_{\mu'_\bullet}) - {\rm dim}({\mathcal Q}_\lambda)] = \langle \rho, |\mu'_\bullet| - \lambda \rangle.
$$

\medskip

When we work in the context of Hecke algebras ${\mathcal H}_q$, the field $k$ will be the finite field $\mathbb F_q$, where $q = p^j$ for a prime $p$.  
In any case, we will always fix a prime $\ell \neq {\rm char}(k)$, and fix an algebraic closure $\overline{\mathbb Q}_\ell$ of $\mathbb Q_\ell$.  We define the dual group $\widehat{G} = \widehat{G}(\overline{\mathbb Q}_\ell)$.  We let $\widehat{T} \subset \widehat{G}$ denote the dual torus of $T$, defined by the equality $X^*(\widehat{T}) = X_*(T)$.

Let $Q^\vee = Q^\vee(G)$ (resp. $Q = Q(G)$) denote the lattice in $X_*(T)$ (resp. $X^*(T)$) spanned by the coroots (resp. roots) of $G$ in $T$.  There is a canonical identification $Q^\vee(G) = Q(\widehat{G})$, by which we can define a notion of simple positive root in $\widehat{G}$ and thus a corresponding Borel subgroup $\widehat{B}$ containing $\widehat{T}$.

When we consider an $r+1$-tuple of coweights $(\mu_\bullet,\lambda)$, it will always be assumed that $\sum_i \mu_i - \lambda \in Q^\vee$.  (When thinking of these as weights of $\widehat{T}$, this amounts to assuming that $\sum_i \mu_i - \lambda \in Q(\widehat{G})$.)

For $\mu$ dominant we let $\Omega(\mu)$ denote the set of weights of the irreducible representation of $\widehat{G}$ with highest weight $\mu$.  For $\nu \in X_*(T)$, we let $S_\nu = Ut_\nu e_0$.  

If $\mu$ is dominant then we denote by $V_\mu$ the irreducible $\widehat{G}$-module with highest weight $\mu$.  Its contragredient $(V_\mu)^*$ is also irreducible, so we can define the {\em dual} dominant coweight $\mu^*$ by the equality $V_{\mu^*} = (V_\mu)^*$.  We have $\mu^* = -w_0 \mu$.

We shall make frequent use of the 
fact that $S_\nu \cap \overline{\mathcal Q}_\mu \neq \emptyset$ only if $\nu \in \Omega(\mu)$ (\cite{BT}, 4.4.4, or \cite{NP}, Lemme 4.2).  For any $\nu \in X_*(T)$, let $\nu_d$ denote the unique $B$-dominant element in $W\nu$.  The Weyl group permutes the set of (co)weights, and we let $W_{\mu}$ denote the stabilizer in $W$ of $\mu$.

\medskip

Recall that a coweight $\mu$ is {\em minuscule} provided that $\langle \alpha, \mu \rangle \in \{-1, 0 ,1 \}$, for every root $\alpha$.  Viewing $\mu$ as a weight of $\widehat{G}$, this is equivalent to the statement that $\Omega(V_\mu) = W\mu$ (see \cite{Bou}).

\subsection{Local triviality of the morphism $m_{\mu_\bullet}$} \label{local_triviality}

Let $X = \cup_{i} X_i$ and $Y = \cup_j Y_j$ be stratifications of algebraic varieties over $k$ by locally closed subvarieties, having the property that the boundary of any stratum is a union of other strata.  

Suppose we have a morphism $f: X \rightarrow Y$.  We suppose that $f$ is proper and that each $f(X_i)$ is a union of strata $Y_j$.  We say $f$ is {\em locally trivial in the stratified sense}, if for every $y \in Y_j$ there is a Zariski-open subset $V \subset Y_j$ with $y \in V$, and a stratified variety $F$, such that there is an isomorphism of stratified varieties
\begin{equation} \label{eq:local_triv}
f^{-1}(V) \cong F \times V
\end{equation}
which commutes with the projections to $V$. \footnote{This definition differs from that used in \cite{H}.  In that paper, a weaker notion of ``locally trivial in the stratified sense'' was used.  This paper requires the present (more conventional) definition.}

The following lemma is well-known, see \cite{MV1}.  We give the proof for the convenience of the reader.

\begin{lemma} \label{local_triv_lemma}
The morphism $m_{\mu_\bullet}$ is Zariski-locally trivial in the stratified sense.
\end{lemma}

\begin{proof}
Fix $y \in {\mathcal Q}_\lambda \subset \overline{\mathcal Q}_{|\mu_\bullet|}$.  We can identify ${\mathcal Q}_\lambda$ with the quotient in the notation of loop groups
$$
{\mathcal Q}_\lambda = L^{\geq 0}G/ L^{\geq 0}G \cap L^{\geq \lambda}G,
$$
where by definition $L^{\geq \lambda}G = \lambda L^{\geq 0}G \lambda^{-1}.$  Suppose that Zariski-locally on the base, the projection
\begin{equation} \label{eq:quotient}
L^{\geq 0}G \rightarrow L^{\geq 0}G / L^{\geq 0}G \cap L^{\geq \lambda}G
\end{equation}
has a section.  Then it is easy to see that Zariski-locally, there is an isomorphism as in 
(\ref{eq:local_triv}) for $f = m_{\mu_\bullet}$.  Indeed, suppose $L_\bullet = (L_1, \dots, L_r) \in \widetilde{\mathcal Q}_{\mu_\bullet}$ has $L_r \in {\mathcal Q}_\lambda$.  Then for $L_r$ in a Zariski-neighborhood $V$ of $e_\lambda$ in ${\mathcal Q}_\lambda$, we can write $L_r = ke_\lambda$ for a well-defined $k \in L^{\geq 0}G$, the image of $L_r$ under the local section.  Then we may define $(\ref{eq:local_triv})$ by
$$
(L_1, \dots, L_r) \mapsto (k^{-1}L_1, \dots, k^{-1}L_r) \times L_r.
$$

It remains to prove that (\ref{eq:quotient}) is Zariski-locally trivial.  By \cite{NP}, Lemme 2.3, we can write 
$L^{\geq 0}G \cap L^{\geq \lambda}G = P_\lambda \ltimes (L^{>0}G \cap L^{\geq \lambda}G)$, where $P_\lambda \subset G$ is the parabolic subgroup corresponding to the roots $\alpha$ such that $\langle \alpha, \lambda \rangle \leq 0$, and where $L^{>0}G$ is the kernel of the morphism $L^{\geq 0}G \rightarrow G$ induced by $t \mapsto 0$.  We also have an obvious isomorphism $L^{\geq 0}G = G \ltimes L^{>0}G$.  Then (\ref{eq:quotient}) can be factored as the composition of two projections
\begin{equation} \label{eq:quotient1}
L^{\geq 0}G \rightarrow L^{\geq 0}G/ L^{>0}G \cap L^{\geq \lambda}G = G \times [L^{>0}G/ L^{>0}G \cap L^{\geq \lambda}G],
\end{equation} 
and 
\begin{equation} \label{eq:quotient2}
G \times [L^{>0}G/ L^{>0}G \cap L^{\geq \lambda}G] \rightarrow L^{\geq 0}G/L^{\geq 0}G \cap L^{\geq \lambda}G.
\end{equation}
Here the first projection is the obvious one, and the second projection is the quotient for the right action of $P_\lambda$ on $G \times [L^{>0}G/L^{>0} \cap L^{\geq \lambda}G]$ given by
$$
(g, g^+ L^{>0}G \cap L^{\geq \lambda}G) \cdot p = (gp, p^{-1}g^+pL^{>0}G \cap L^{\geq \lambda}G).
$$

The morphism (\ref{eq:quotient1}) is actually trivial, because the multiplication map
$$
(L^{>0}G \cap L^{< \lambda}G)  \times (L^{>0}G \cap L^{\geq \lambda}G) \rightarrow L^{>0}G
$$
is an isomorphism, where $L^{<\lambda}G := \lambda L^{<0}G \lambda^{-1}$ and where $L^{<0}G$ denotes the kernel of the map $G(k[t^{-1}]) \rightarrow G$ induced by $t^{-1} \mapsto 0$; see \cite{NP}, $\S 2$.

The morphism (\ref{eq:quotient2}) has local sections in the Zariski topology, coming from the embedding of the ``big cell'' $U_{\overline{P_\lambda}} \hookrightarrow G/P_\lambda$.  This completes the proof.
\end{proof}

\subsection{Review of information carried by fibers of convolution morphisms} \label{content_of_fibers}

The following well-known result plays a key role in this article.  

\begin{theorem}[Geometric Satake Isomorphism -- weak form] \label{weak_GS}
For every tuple $(\mu_\bullet,\lambda)$, and every $y \in {\mathcal Q}_\lambda$, there is an equality 
$$
{\rm dim}(V^\lambda_{\mu_\bullet}) = \# \, \mbox{irreducible components of $m^{-1}_{\mu_\bullet}(y)$ having dimension $\langle \rho, |\mu_\bullet| -\lambda \rangle$.}
$$
\end{theorem}

See \cite{H}, $\S 3$ for the proof of this assuming the geometric Satake isomorphism in the context of finite residue fields.  We also have the following elementary lemma.

\begin{lemma} \label{c(q)}
For $(\mu_\bullet, \lambda), y$ as above,
$$
c^\lambda_{\mu_\bullet}(q) = \# \, ({\mathcal Q}_{\mu_\bullet} \cap m^{-1}_{\mu_\bullet}(y))({\mathbb F}_q).
$$
\end{lemma}

For context we recall following \cite{H} that the above two statements together with the Weil conjectures yield the following expression for the Hecke algebra structure constants.

\begin{cor} [\cite{KLM}] \label{structure_constants_formula}  With $\mu_\bullet, \lambda$ as above, the Hecke algebra structure constant is given by the formula
$$
c^\lambda_{\mu_\bullet}(q) = {\rm dim}(V^\lambda_{\mu_\bullet}) \, q^{\langle \rho, |\mu_\bullet|-\lambda \rangle} + (\mbox{terms with lower $q$-degree}).
$$
\end{cor}

This formula was first proved by Kapovich, Leeb, and Millson \cite{KLM}, who deduced it from the results of Lusztig \cite{Lu2}.  It actually provides an algorithm to compute the multiplicities ${\rm dim}(V^\lambda_{\mu_\bullet})$.  Indeed, one can determine the polynomial $c^\lambda_{\mu_\bullet}(q)$ by computing products in an Iwahori-Hecke algebra, using the Iwahori-Matsumoto presentation of that algebra.  Of course this involves the computation of much more than just the leading term of $c^\lambda_{\mu_\bullet}$, so in practice this procedure is not a very efficient way to compute ${\rm dim}(V^\lambda_{\mu_\bullet})$.

However the formula does make it clear that the dimensionality of the fiber $m_{\mu_\bullet}^{-1}(y)$ plays a role in linking the non-vanishing of the structure constants:  if ${\rm dim}(V^\lambda_{\mu_\bullet}) > 0$, then evidently $c^\lambda_{\mu_\bullet}(q) \neq 0$ for all large $q$ (and thus all $q$, by the argument in \cite{H}, $\S 4$).  On the other hand, if $c^\lambda_{\mu_\bullet}(q) \neq 0$, it could well happen that the leading coefficient ${\rm dim}(V^\lambda_{\mu_\bullet})$ is zero.  However, if we knew {\em a priori} that whenever the space 
$m_{\mu_\bullet}^{-1}(y) \cap {\mathcal Q}_{\mu_\bullet}$ is non-empty, it is actually of dimension $\langle \rho, |\mu_\bullet| - \lambda \rangle$, then the non-vanishing of $c^\lambda_{\mu_\bullet}(q)$ would imply the non-vanishing of its leading coefficient.  We will prove this dimension statement for $m_{\mu_\bullet}^{-1}(y) \cap {\mathcal Q}_{\mu_\bullet}$ in the case where each $\mu_i$ is a sum of minuscules, by a reduction to the case where each $\mu_i$ is minuscule.  But as is seen in the reduction step (the ``pulling apart'' Lemma \ref{reduction_to_sums}) it is necessary to prove the stronger fact that in that case, the fibers are not just of largest possible dimension, but are also {\em equidimensional}.  

Our first goal, therefore, is to establish the (equi)dimensionality statement just mentioned (in Theorem \ref{equidim} below).  Let us first pause to mention some related work in the literature.  After the seminal work of Mirkovic-Vilonen \cite{MV1,MV2} on which everything else is based, the author was particularly inspired by the work of Ng\^{o}-Polo \cite{NP}.  Many other authors have had the idea to use the fibers of the morphisms $m_{\mu_\bullet}$ to derive representation-theoretic consequences, and the works of Gaussent-Littelmann \cite{GL} and of J. Anderson \cite{A} seem particularly related to the present one.  In fact, in \cite{A} Anderson independently observed the relation between fibers of convolution morphisms and MV cycles (loc. cit. Theorem 8), which was the starting point in the proof of our Theorem \ref{equidim}.

\section{Equidimensionality of minuscule convolutions} 

\subsection{Proof of the main theorem}
For this section we fix an $r$-tuple $\mu_\bullet = (\mu_1, \dots, \mu_r)$ such that each $\mu_i$ is dominant and {\em minuscule}.  The main result of this paper is the following theorem.

\begin{theorem} \label{equidim}
The fibers of the morphism
$$
m_{\mu_\bullet}: \widetilde{\mathcal Q}_{\mu_\bullet} \rightarrow 
\overline{\mathcal Q}_{|\mu_\bullet|}
$$
are equidimensional.  More precisely, if $y \in {\mathcal Q}_\lambda \subset 
\overline{\mathcal Q}_{|\mu_\bullet|}$, then every irreducible component of 
$m^{-1}_{\mu_\bullet}(y)$ has dimension $\langle \rho, |\mu_\bullet| -\lambda \rangle$.
\end{theorem}

\bigskip

We will prove the theorem by induction on $r$, the number of elements in the tuple $\mu_\bullet$ (the case of $r=1$ being trivial).  Let us suppose the theorem is true for the morphism 
$$
m_{\mu'_\bullet} : {\mathcal Q}_{\mu_1} \tilde{\times} 
\cdots \tilde{\times} {\mathcal Q}_{\mu_{r-1}} \rightarrow 
\overline{\mathcal Q}_{|\mu'_\bullet|}
$$
attached to the $(r-1)$-tuple $\mu'_\bullet = (\mu_1, \dots, \mu_{r-1})$.  We fix an orbit ${\mathcal Q}_\lambda \subset \overline{\mathcal Q}_{|\mu_\bullet|}$ and a point $y \in {\mathcal Q}_\lambda$; we want to prove the equidimensionality of the variety $m_{\mu_\bullet}^{-1}(y)$.  By equivariance under the $K$-action, we can assume that $y = e_\lambda := t_\lambda e_0$.
We suppose that $(L_1, \dots, L_{r-1},L_r) \in m_{\mu_\bullet}^{-1}(y)$.  
Since the relative position ${\rm inv}(L_{r-1},L_r) = \mu_r$ and $L_r = e_\lambda$, we deduce that $L_{r-1} \in t_\lambda K t_\mu e_0$, where $\mu:= \mu^*_r$.  Thus, $L_{r-1}$ ranges over the set 
$$\overline{\mathcal Q}_{|\mu'_\bullet|} \cap 
t_\lambda {\mathcal Q}_\mu.
$$  
Most of the work in the proof of Theorem \ref{equidim} involves the attempt to understand this set.  It is hard to understand the whole set, but as we shall see below, we can exhaust it by locally closed subsets $Z_w$ which are easier to understand.  The locally closed subsets help us compute dimensions of components in $m_{\mu_\bullet}^{-1}(y)$ because, as we shall see, 
\begin{itemize}
\item the morphism $m_{\mu'_\bullet}$ becomes trivial over each of these subsets, and 
\item we can explicitly calculate the dimensions of the subsets.
\end{itemize}

More precisely, since $\mu$ is a minuscule coweight, we decompose $t_\lambda {\mathcal Q}_\mu$ as the union of the locally closed subsets 
$$
t_\lambda(\, U t_{w\mu}e_0 \cap K t_{\mu}e_0) = \,  U t_{\lambda + w\mu} e_0 \cap t_\lambda K t_{\mu} e_0,
$$
where $w \in W/W_\mu$.  
\medskip

Using the above decomposition of $t_\lambda {\mathcal Q}_\mu$, we see that $\overline{\mathcal Q}_{|\mu'_\bullet|} \cap t_\lambda{\mathcal Q}_\mu$ is the disjoint union of the following locally-closed subvarieties
$$
Z_w := \overline{\mathcal Q}_{|\mu'_\bullet|} \cap t_\lambda(S_{w\mu} \cap {\mathcal Q}_\mu),
$$
where we recall that $S_\nu := Ut_\nu e_0$ for $\nu \in X_*(T)$.  For brevity, we let $X_w := S_{w\mu} \cap {\mathcal Q}_\mu$, so that $Z_w = \overline{\mathcal Q}_{|\mu'_\bullet|} \cap t_\lambda X_w$.

\medskip

We begin the proof with a preliminary lemma concerning $X_w$, and then proceed to some lemmas concerning $Z_w$. 

The following result is due to Ng\^{o}-Polo \cite{NP}.  In this statement $R^+$ denotes the set of positive roots relative to the Borel subgroup $B = TU$.

\begin{lemma} \label{intersect}
Let $\nu$ be and dominant coweight, and let $w \in W$.  Then there is an isomorphism of varieties 
$$
S_{w\nu} \cap {\mathcal Q}_\nu = \prod_{\alpha \in R^+ \cap wR^+} \prod_{i=0}^{\langle w^{-1}\alpha, \nu \rangle -1} U_{\alpha,i} t_{w\nu}e_0,
$$ 
where $U_\alpha \subset LU$ is the root subgroup corresponding to the positive root $\alpha$, and $U_{\alpha,i}$ consists of the elements in $U_{\alpha}$ of form $u_\alpha(x t^i)$, $x \in k$, where $u_\alpha : {\mathbb G}_a \rightarrow U_\alpha$ is the root homomorphism for $\alpha$.  

In particular, applying this to $\nu = \mu$, we see that $X_w$ is an affine space of dimension
$${\rm dim}(X_w) = \langle \rho, \mu + w \mu \rangle.$$
\end{lemma}

\begin{proof}
The isomorphism is proved in \cite{NP}, Lemme 5.2.  The dimension formula then follows, using the 
formula 
$$
\rho + w^{-1}\rho = \sum_{\alpha \in R^+ \cap w^{-1}R^+} \, \alpha,
$$
and its consequence
$$
\langle \rho, \mu + w\mu \rangle = \sum_{\alpha \in R^+ \cap wR^+} \langle \alpha, w\mu \rangle.
$$
\end{proof}

\medskip

Consider next the morphism 
$$
p: m_{\mu_\bullet}^{-1}(t_\lambda e_0) \twoheadrightarrow \overline{\mathcal Q}_{|\mu'_\bullet|} \cap t_\lambda {\mathcal Q}_{\mu}
$$
given by $p(L_1,\dots, L_{r-1},L_r) = L_{r-1}$.

\begin{lemma} \label{lemma1}  Suppose $Z_w \neq \emptyset$.  Then $t_{\lambda + w\mu}e_0 \in Z_w$, and 
the morphism $p : p^{-1}(Z_w) \twoheadrightarrow Z_w$ is trivial.  In particular 
\begin{equation} \label{eq:triviality_formula}
p^{-1}(Z_w) \cong m_{\mu'_\bullet}^{-1}(t_{\lambda + w\mu} e_0) \times Z_w.
\end{equation}
\end{lemma}

\begin{proof}  
Note that 
\begin{align*}
\overline{\mathcal Q}_{|\mu'_\bullet|} \cap t_{\lambda}(Ut_{w\mu}e_0 \cap Kt_\mu e_0) \neq \emptyset &\Rightarrow \overline{\mathcal Q}_{|\mu'_\bullet|} \cap Ut_{\lambda + w\mu} e_0 \neq \emptyset \\
&\Rightarrow \lambda + w\mu \in \Omega(|\mu'_\bullet|) \\
&\Rightarrow t_{\lambda + w\mu}e_0 \in \overline{\mathcal Q}_{|\mu'_\bullet|}.
\end{align*}
Since clearly $t_{\lambda + w\mu}e_0 \in t_\lambda(S_{w\mu} \cap {\mathcal Q}_\mu)$, 
the first statement follows.

Now we prove the second statement.  By Lemma \ref{intersect}, an element in $t_\lambda(Ut_{w\mu} e_0 \cap Kt_{\mu} e_0)$ can be written uniquely in the form 
\begin{equation} \label{eq:t_lambda_X_w}
t_\lambda( \prod_{\stackrel{\alpha \in R^+}{\langle \alpha, w\mu \rangle = 1}} u_\alpha(x_\alpha)t_{w\mu} e_0),
\end{equation}
where $x_\alpha \in k$ for all $\alpha$.  Because $\lambda$ is dominant, we can write this as 
\begin{equation} \label{eq:n_0}
n_0 t_{\lambda + w\mu}e_0,
\end{equation}
for some uniquely determined $n_0 \in L^{\geq 0}U$.   Therefore if $(L_1, \dots, L_{r-1},L_r) \in p^{-1}(Z_w)$, we can write $L_{r-1} = n_0 t_{\lambda + w\mu} e_0$ for a uniquely determined $n_0 \in L^{\geq 0}U$.  Then the isomorphism 
$p^{-1}(Z_w) \rightarrow m_{\mu'_\bullet}^{-1}(t_{\lambda + w\mu} e_0) \times Z_w$ is the map sending
$$
(L_1, \dots, L_{r-1},L_r) \mapsto (n^{-1}_0L_1, \dots, n^{-1}_0L_{r-1}) \times L_{r-1}.
$$ 

Note that since $m_{\mu'_\bullet}$ is $K$-equivariant and $n_0 \in K$, we also have an isomorphism
$$
p^{-1}(Z_w) \cong m_{\mu'_\bullet}^{-1}(L_{r-1}) \times Z_w,
$$
for any $L_{r-1} \in Z_w$.
\end{proof}

\begin{lemma}\label{when_Z_w_is_non-empty}
\begin{enumerate}
\item We have $t_\lambda X_w \subset \mathcal Q_{(\lambda + w\mu)_d}$.  
\item We have $Z_w \neq \emptyset \Longleftrightarrow (\lambda + w\mu)_d \leq |\mu'_\bullet|$, in which case $Z_w = t_\lambda X_w$.
\end{enumerate}
\end{lemma}

\begin{proof}
\noindent (1):  This follows immediately from equations (\ref{eq:t_lambda_X_w}) and (\ref{eq:n_0}) above.

\noindent (2):  Clearly since $t_\lambda X_w \subset \mathcal Q_{(\lambda + w\mu)_d}$, we have $Z_w \neq \emptyset$ if and only if $\mathcal Q_{(\lambda + w\mu)_d} \subset \overline{\mathcal Q}_{|\mu'_\bullet|}$, which holds if and only if $(\lambda + w\mu)_d \leq |\mu'_\bullet|.$  It is also clear that $Z_w = t_\lambda X_w$ in that case.
\end{proof}

\medskip

Recall that $p^{-1}(Z_w) \subset m_{\mu_\bullet}^{-1}(e_\lambda)$.  By the semi-smallness of $m_{\mu_\bullet}$, for every $w$ we have ${\rm dim}(p^{-1}(Z_w)) \leq \langle \rho, |\mu_\bullet| - \lambda \rangle$.  We call $Z_w$ {\em good} if equality holds.  If $Z_w$ is good, then $p^{-1}(Z_w)$ is equidimensional of dimension $\langle \rho, |\mu_\bullet| - \lambda \rangle$.  If $Z_w$ is not good, then $p^{-1}(Z_w)$ is a equidimensional of strictly smaller dimension.  (We already know using the induction hypothesis and triviality that $p^{-1}(Z_w)$ is equidimensional for every $w$.)   

We need to give a concrete criterion for ``$Z_w$ is good''.  

\begin{lemma} \label{when_Z_w_is_good}
Suppose $Z_w \neq \emptyset$.  Then $Z_w$ is good if and only if $\lambda + w\mu$ is dominant.
\end{lemma}

\begin{proof}
Using (\ref{eq:triviality_formula}) we have 
$$
{\rm dim}(Z_w) + {\rm dim}(m_{\mu'_\bullet}^{-1}(t_{(\lambda + w\mu)_d}e_0)) = {\rm dim}(p^{-1}(Z_w)) \leq \langle \rho, |\mu_\bullet| - \lambda \rangle,
$$
with equality if and only if $Z_w$ is good.  Since $Z_w$ is non-empty, we know by Lemma \ref{when_Z_w_is_non-empty} and Lemma \ref{intersect} that ${\rm dim}(Z_w) = 
\langle \rho, \mu + w\mu \rangle$.  Using this together with our induction hypothesis that ${\rm dim}(m_{\mu'_\bullet}^{-1}(t_{(\lambda + w\mu)_d}e_0)) = \langle \rho, |\mu'_\bullet| - (\lambda + w\mu)_d \rangle$, and the equality $\langle \rho, |\mu_\bullet|  \rangle = \langle \rho, |\mu'_\bullet| + \mu \rangle$, the above statement becomes
$$
\langle \rho, (\lambda + w\mu) - (\lambda + w\mu)_d \rangle \leq 0
$$
with equality if and only if $Z_w$ is good.

Thus we see that $\lambda + w\mu = (\lambda + w\mu)_d$ if and only if $Z_w$ is good.
\end{proof}

\begin{lemma} \label{first_closure_relation}
For any $w \in W$, let $\overline{Z_w}$ denote the closure of $Z_w$ in $\overline{\mathcal Q}_{|\mu'_\bullet|} \cap t_\lambda{\mathcal Q}_\mu$.  Suppose $Z_w$ and $Z_{w'}$ are non-empty.  Then $Z_w \subset \overline{Z_{w'}}$ if $w' \leq w$ in the Bruhat order on $W$.
\end{lemma}

\begin{proof}
Let $P^-$ denote the standard parabolic determined by the set of roots satisfying $\langle \alpha, \mu \rangle \leq 0$.  By Lemma \ref{when_Z_w_is_non-empty}, we have $Z_{w'} = t_\lambda X_{w'}$ and $Z_w = t_\lambda X_w$, so that the closure relations for the $Z_w$'s inside $t_\lambda {\mathcal Q}_\mu$ are determined by those for the $X_w$'s inside ${\mathcal Q}_\mu$.  

By Lemma \ref{intersect} the ``reduction modulo $t$'' isomorphism
${\mathcal Q}_\mu ~\widetilde{\rightarrow} ~ G/P^-$ induces an isomorphism
$$
X_w ~ \widetilde{\rightarrow} ~ UwP^-/P^-.
$$
(comp. \cite{NP}, Lemme 6.2.)  
The result now follows from the relation between the Bruhat order on $W$ and the closure relations for $U$-orbits in $G/P^-$.
\end{proof}

\medskip

\noindent {\em Proof of Theorem \ref{equidim}}.  Consider again the morphism
$$
p: m_{\mu_\bullet}^{-1}(t_\lambda e_0) \rightarrow \overline{\mathcal Q}_{|\mu'_\bullet|} \cap t_\lambda {\mathcal Q}_{\mu}
$$
given by $p(L_1,\dots, L_{r-1},L_r) = L_{r-1}$. 
We have constructed a decomposition of the range by locally closed sets 
$Z_w$, $w \in W$, over which $p$ is trivial.  Some of the sets $Z_w$ might be empty, but for non-empty $Z_w$, we now have a useful description of those which are good (Lemma \ref{when_Z_w_is_good}).  A priori we do not know whether any good subsets exist, but in the course of the proof we shall see that they do.

Using our induction hypothesis, we know that for good $Z_w$, the set $\overline{p^{-1}(Z_w)}$ is a union of irreducible components of $m_{\mu_\bullet}^{-1}(t_\lambda e_0)$ having dimension $\langle \rho, |\mu_\bullet| - \lambda \rangle$.  It remains to prove that for any non-empty set $Z_w$ which is {\em not} good, there exists a non-empty good $Z_{w^*}$ such that $p^{-1}(Z_w) \subset \overline{p^{-1}(Z_{w^*})}$.

The first step is to find a good $Z_{w^*}$ such that $Z_w \subset \overline{Z_{w^*}}$.  Assume $Z_w$ is non-empty but is {\em not} good.  Then let $w^*$ be the unique element of minimal length in the subset $W_\lambda w W_\mu$ of $W$.  Since $\lambda + w^* \mu$ and $\lambda + w\mu$ are $W_\lambda$-conjugate, we have $(\lambda + w^*\mu)_d = (\lambda + w\mu)_d$ and hence by Lemma \ref{when_Z_w_is_non-empty}, $Z_{w^*} \neq \emptyset$.  Then by Lemma \ref{first_closure_relation} we have $Z_w \subset \overline{Z_{w^*}}$.  It remains to prove $Z_{w^*}$ is good, i.e. that $\lambda + w^*\mu$ is dominant (Lemma \ref{when_Z_w_is_good}).  But if $\lambda + w^*\mu$ is not dominant, there is a positive root $\alpha$ with $\langle \alpha, \lambda + w^*\mu \rangle <0$.  Since $\lambda$ is dominant and $w^*\mu$ is minuscule, we must have $\langle \alpha, \lambda \rangle = 0$ and $\langle \alpha, w^*\mu \rangle = -1$.  The latter implies that $(w^*)^{-1}\alpha < 0$, which means that $s_\alpha w^* < w^*$ in the Bruhat order on $W$.  But since $s_\alpha \in W_\lambda$, this contradicts the definition of $w^*$.

To complete the proof of Theorem \ref{equidim}, we need to show that $p^{-1}(Z_w) \subset \overline{p^{-1}(Z_{w^*})}$.  Roughly, this follows because $\lambda + w\mu$ and $\lambda + w^*\mu$ are $W_\lambda$-conjugate, hence both $Z_w$ and $Z_{w^*}$ belong to ${\mathcal Q}_{\lambda + w^*\mu}$, over which $m_{\mu'_\bullet}$ is locally trivial (Lemma \ref{local_triv_lemma}).  
More precisely, suppose 
$L_\bullet := (L_1, \dots, L_{r-1}, L_r) \in p^{-1}(Z_w)$.  Write $y \in Z_w$ for $L_{r-1}$, the image of $(L_1, \dots, L_{r-1})$ under $m_{\mu'_\bullet}$.  Let $F = m^{-1}_{\mu'_\bullet}(y)$.  Choose an open neighborhood $y \in U \subset {\mathcal Q}_{\lambda + w^*\mu} \cap t_\lambda{\mathcal Q}_\mu$ over which $m_{\mu'_\bullet}$ is trivial.  Since $p$ is the just the restriction of $m_{\mu'_\bullet}$ over 
$\overline{\mathcal Q}_{|\mu'_\bullet|} \cap t_\lambda {\mathcal Q}_\mu$, it follows that $p$ is also trivial over $U$, so that $p^{-1}(U) \cong F \times U$.  

To show 
$$L_\bullet \in \overline{p^{-1}(Z_{w^*})},
$$
 it is enough to show
$$L_\bullet \in p^{-1}(U) \cap \overline{p^{-1}(Z_{w^*})}.
$$
The intersection on the right hand side contains
\begin{align*}
\overline{p^{-1}(U) \cap p^{-1}(Z_{w^*})} &\cong \overline{F \times (U \cap Z_{w^*})} \\
&= F \times \overline{U \cap Z_{w^*}} \\
&= F \times U \\
&\cong p^{-1}(U),
\end{align*}
where for $V$ open and $A$ arbitrary, $\overline{V \cap A}$ denotes the closure of $V \cap A$ in the subspace topology on $V$.  In proving the second equality we have used the fact that $U \cap Z_{w^*}$ is non-empty and open in $U$, and that $U$ is irreducible.  These statements follow from the fact that the irreducible set $Z_{w^*}$ is open and dense in $Z_{[w^*]} := {\mathcal Q}_{\lambda + w^*\mu} \cap t_\lambda{\mathcal Q}_\mu$ (as proved in Lemma \ref{lemma7} below).

Our assertion now follows since $L_\bullet$ obviously belongs to $p^{-1}(U)$.  This completes the proof of Theorem \ref{equidim}, modulo  Lemma \ref{lemma7} below.  
\qed

\bigskip

\subsection{Description of closure relations, and paving by affine spaces} \label{paving_by_affines_section}

Lemma \ref{first_closure_relation} gives a partial description of the closure relations between the $Z_w$ subsets. 
Our present aim is to give a complete description.  

Note that every class $[w] \in W_\lambda \backslash W/W_\mu$ gives rise to a well-defined $K$-orbit ${\mathcal Q}_{\lambda + w\mu}$.  Each double coset is represented by a unique element $w^*$ of minimal length.  In other words, $w^*$ is the unique element of minimal length in its double coset $W_\lambda w^* W_\mu$.  As remarked in the proof of Theorem \ref{equidim}, the coweight $\lambda + w^*\mu$ is dominant.

If $\lambda + w^* \mu \leq |\mu'_\bullet|$, we denote 
$$
Z_{[w^*]} := {\mathcal Q}_{\lambda + w^*\mu} \cap t_\lambda{\mathcal Q}_\mu = \bigcup_{w \in W_\lambda w^* W_\mu} {\mathcal Q}_{\lambda + w^*\mu} \cap t_\lambda(S_{w\mu} 
\cap {\mathcal Q}_\mu);
$$
in case $\lambda + w^*\mu \nleq |\mu'_\bullet|$, set $Z_{[w^*]} = \emptyset$.
Concerning the second equality defining $Z_{[w^*]}$, it is clear that the left hand side contains the right hand side.  To prove the other inclusion, note that if a subset of the left hand side of form ${\mathcal Q}_{\lambda + w^*\mu} \cap t_\lambda(S_{w\mu} \cap {\mathcal Q}_\mu)$ is non-empty, then it is a $Z_w$, and has $(\lambda + w\mu)_d = \lambda + w^*\mu$, from which it follows that $w \in W_\lambda  w^* W_\mu$.

Clearly we have a decomposition by locally closed (possibly empty) subsets 
$$
\overline{\mathcal Q}_{|\mu'_\bullet|} \cap t_\lambda {\mathcal Q}_\mu = \coprod_{[w^*] \in W_\lambda \backslash W/W_\mu} Z_{[w^*]}.
$$

\begin{lemma} \label{lemma7}
The following statements hold.
\begin{enumerate}
\item[(a)]  For $w \in W$, let $\overline{Z_w}$ denote the closure of $Z_w$ in $\overline{\mathcal Q}_{|\mu'_\bullet|} \cap t_\lambda{\mathcal Q}_\mu$.  If $Z_w \neq \emptyset$, then 
$$
\overline{Z_w} = \bigcup_{v \geq w} Z_v.
$$
Furthermore, $Z_w \neq \emptyset \Rightarrow Z_v \neq \emptyset$ for all $v \geq w$.
\item[(b)]  We have $Z_{w^*} \neq \emptyset$ if and only if $Z_w \neq \emptyset$ for any $w \in [w^*]$.
\item[(c)]  If $Z_{w^*} \neq \emptyset$, then the map $t_{\lambda}{\mathcal Q}_\mu \rightarrow {\mathcal Q}_\mu \rightarrow G/P^-$ induces an isomorphism
\begin{equation*}
Z_{[w^*]} ~ \widetilde{\longrightarrow} ~ \bigcup_{w \in W_\lambda w^*W_\mu} U w P^-/P^-,
\end{equation*}
and furthermore $Z_{w^*}$ is dense and open in $Z_{[w^*]}$. 
\item[(d)] Let $v^*$ denote a minimal representative for a double coset $W_\lambda v^* W_\mu$.  
If $Z_{w} \neq \emptyset$ then 
$$
\overline{Z_w} \cap Z_{[v^*]} = \bigcup_{v \in [v^*], \,\, v \geq w} Z_v.
$$
\end{enumerate}
\end{lemma}

\begin{cor} \label{components}
The irreducible components of $\overline{\mathcal Q}_{|\mu'_\bullet|} \cap t_\lambda{\mathcal Q}_\mu$ are the closures $\overline{Z_{w^*}}$, where $w^*$ ranges over the minimal elements in the set $\{ v^*  ~|~ \lambda + v^*\mu \leq |\mu'_\bullet| \}$.
\end{cor}

\begin{proof}

\noindent
(a): The morphism  
$$
\overline{\mathcal Q}_{|\mu'_\bullet|} \cap t_\lambda{\mathcal Q}_\mu \rightarrow G/P^-
$$
is a closed immersion, hence proper.  So if $Z_w \neq \emptyset$, then the image of $\overline{Z_w}$ is the closure 
$$
\overline{UwP^-/P^-} = \bigcup_{v \geq w} UvP^-/P^-.
$$
It follows that $Z_w \neq \emptyset \Rightarrow Z_v \neq \emptyset$, for all $v \geq w$, and that the closure 
above is the image of 
$$
\bigcup_{v \geq w} Z_v.
$$

 (b):  This follows from Lemma \ref{when_Z_w_is_non-empty}, using the equality $(\lambda + w\mu)_d = \lambda + w^*\mu$.

(c):  This is easy, the main point being that $Uw^*P^-/P^-$ is clearly open and dense in the union of {\em all} $UwP^-/P^-$ for $w \in W$ with $w \geq w^*$.    

(d): This follows from (a)-(c).
\end{proof}

\noindent {\em Proof of Corollary \ref{paved_by_affines}.}
We can now prove that $m_{\mu_\bullet}^{-1}(y)$ is indeed paved by affine spaces.  Let us recall what this means.  By definition, a scheme $X$ is {\em paved by affine spaces} if there is an increasing filtration $\emptyset = X_0 \subset X_1 \subset \cdots \subset X_n = X$ by closed subschemes $X_i$ such that each successive difference $X_i \backslash X_{i-1}$ is a (topological) disjoint union of affine spaces 
${\mathbb A}^{n_{ij}}$. 

We prove the corollary by induction on $r$: assume every fiber of $m_{\mu'_\bullet}$ is paved by affines.  
By the above discussion, $\overline{\mathcal Q}_{|\mu'_\bullet|} \cap t_\lambda{\mathcal Q}_\mu$ is a disjoint union of certain (non-empty) locally closed subsets $Z_w$, each of which is isomorphic to an affine 
space.  The boundary of each such $Z_w$ is a union of other strata $Z_v$.  The triviality statement of Lemma \ref{lemma1} and the induction hypothesis then shows that each variety $p^{-1}(Z_w)$ is paved by affine spaces.   

These remarks imply (by an inductive argument) that $m_{\mu_\bullet}^{-1}(y)$ is paved by affine spaces. \qed

\medskip

\begin{Question} \label{questions}  Suppose $\mu_\bullet$ is a general $r$-tuple of coweights $\mu_i$ (not necessarily minuscule).  Which fibers $m_{\mu_\bullet}^{-1}(y)$  are paved by affine spaces?  Does every fiber $m_{\mu_\bullet}^{-1}(y)$ admit a Hessenberg paving, in the sense of \cite{GKM}, $\S 1$?
\end{Question}

\medskip

\section{Equidimensionality results for sums of minuscules} \label{further_equidim_section}

This section concerns what we can say when the $\mu_i$'s are not all minuscule.  We will consider the fibers of $m_{\mu_\bullet}$, where each $\mu_i$ is a sum of dominant minuscule coweights.  Assume $y \in {\mathcal Q}_\lambda \subset \overline{\mathcal Q}_{|\mu_\bullet|}$.

\begin{prop} \label{further_equidim}  
\begin{enumerate} 
\item [(1)]  Let ${\mathcal Q}_{\mu'_\bullet} \subset \widetilde{\mathcal Q}_{\mu_\bullet}$ be the stratum indexed by $\mu'_\bullet = (\mu'_1, \dots, \mu'_r)$ for dominant coweights $\mu'_i \preceq \mu_i$ ($1 \leq i \leq r$).  Then any irreducible component of the fiber $m_{\mu_\bullet}^{-1}(y)$ whose generic point belongs to ${\mathcal Q}_{\mu'_\bullet}$ has dimension $\langle \rho, |\mu'_\bullet| - \lambda \rangle$.  
\item [(2)]  Suppose each $\mu'_i$ is a sum of dominant minuscule coweights.  Suppose that  
$$
m_{\mu_\bullet}^{-1}(y) \cap {\mathcal Q}_{\mu'_\bullet} = m_{\mu'_\bullet}^{-1}(y) \cap {\mathcal Q}_{\mu'_\bullet}
$$
is non-empty.  Then $m_{\mu_\bullet}^{-1}(y) \cap {\mathcal Q}_{\mu'_\bullet}$ is equidimensional of dimension
$$
{\rm dim}(m_{\mu_\bullet}^{-1}(y) \cap {\mathcal Q}_{\mu'_\bullet}) = \langle \rho, |\mu'_\bullet| - \lambda \rangle.$$
\end{enumerate}
\end{prop}

Note that Theorem \ref{thm_B} follows from part (2), if we take $\mu'_\bullet = \mu_\bullet$.

\begin{proof}
Part (2) follows from part (1).  Part (1) follows from Theorem \ref{equidim} and the following lemma, whose proof appears in \cite{H} (Proof of Prop. 1.8).
\end{proof}

\begin{lemma}[The pulling apart lemma]  \label{reduction_to_sums}
Suppose $\mu_i = \sum_j \nu_{ij}$, for each $i$, and consider the diagram 
$$
\xymatrix{
\widetilde{\mathcal Q}_{\nu_{\bullet \bullet}} \ar[r]^{\eta} & \widetilde{\mathcal Q}_{\mu_\bullet} \ar[r]^{m_{\mu_\bullet}} & \overline{\mathcal Q}_{|\mu_\bullet|},
}
$$
where $\eta = m_{\nu_{1 \bullet}} \tilde{\times} \cdots \tilde{\times} m_{\nu_{r \bullet}}$ and hence $m_{\mu_\bullet} \circ \eta = m_{\nu_{\bullet \bullet}}$.  Then if $m_{\nu_{\bullet \bullet}}^{-1}(e_\lambda)$ is {\em equidimensional} of dimension $\langle \rho, |\nu_{\bullet \bullet}| - \lambda \rangle$, the morphism $m_{\mu_{\bullet}}$ satisfies the conclusion of Proposition \ref{further_equidim}, part (1).
\end{lemma}

\smallskip

\begin{Remark} \label{quasi_min_remark}  In general, the fiber $m^{-1}_{\mu_\bullet}(e_\lambda)$ is not equidimensional of dimension $\langle \rho, |\mu_\bullet| - \lambda \rangle$.  Following \cite{KLM}, $\S 9.5$, consider for example the group $G = {\rm SO}_5$ (so $\widehat{G} = {\rm Sp}_4(\mathbb C)$), where one fundamental weight of $\widehat{G}$ is minuscule and the other is quasi-minuscule ($\mu$ is {\em quasi-minuscule} if $\Omega(V_{\mu}) = W\mu \cup \{0\}$).  The implication ${\rm Hecke}(\mu_\bullet,\lambda) \Rightarrow {\rm Rep}(\mu_\bullet,\lambda)$ does not always hold.  In fact, let 
$$
\mu_1 = \mu_2 = \mu_3 = \alpha_1 + \alpha_2 = (1,1),
$$ where $\alpha_i$ are the two simple roots of $\widehat{G}$, following the conventions of \cite{Bou}.  Let $\lambda = 0$.  In \cite{KLM} it is shown that $V^\lambda_{\mu_\bullet} = 0$ and $c^\lambda_{\mu_\bullet} = q^5 - q \neq 0$. 

We see using Lemma \ref{c(q)} that 
$$
{\rm dim}(m_{\mu_\bullet}^{-1}(e_0)) = {\rm dim}(m_{\mu_\bullet}^{-1}(e_0) \cap {\mathcal Q}_{\mu_\bullet}) = 5
$$
which is strictly less than $\langle \rho, |\mu_\bullet| \rangle = 6$.

Since every coweight of ${\rm SO}_5$ is a sum of minuscule and quasi-minuscule coweights, this example together with Lemma \ref{reduction_to_sums} yields: if we assume each $\mu_i$ is minuscule or quasi-minuscule, in general the fibers $m_{\mu_\bullet}^{-1}(y)$ are {\em not} all equidimensional of the maximal possible dimension.
\end{Remark}

\medskip

\section{Relating structure constants for sums of minuscules}

\begin{cor} \label{structure_constants}
If every $\mu_i$ is a sum of minuscules, then
$$
{\rm Rep}(\mu_\bullet, \lambda) \Leftrightarrow {\rm Hecke}(\mu_\bullet, \lambda).
$$
\end{cor}

\begin{proof}
The argument is as in \cite{H}, which handled the case of ${\rm GL}_n$.  Namely, we prove the implication $\Leftarrow$ as follows.  If ${\rm Hecke}(\mu_\bullet, \lambda)$ holds, then $m_{\mu_\bullet}^{-1}(e_\lambda) \cap {\mathcal Q}_{\mu_\bullet} \neq \emptyset$, and then by Proposition \ref{further_equidim} (2), we see that the dimension of this intersection is $\langle \rho, |\mu_\bullet| - \lambda \rangle$.  Hence by Theorem \ref{weak_GS}, the property ${\rm Rep}(\mu_\bullet, \lambda)$ holds.  
\end{proof}

\begin{Remark} \label{no_assumption_on_lambda}
Note that there is no assumption on the coweight $\lambda$.  In particular, $\lambda$ need not be a sum of dominant minuscule coweights.  After this result was obtained, an improvement was announced in a preprint of Kapovich-Millson \cite{KM}, for the case $r=2$.  This improvement states that 
$$
{\rm Hecke}(\mu_1, \mu_2, \lambda) \Rightarrow {\rm Rep}(\mu_1, \mu_2, \lambda)
$$
as long as at least {\em one} of the coweights $\mu_1, \mu_2$ or $\lambda$ is a sum of minuscules (instead of two of them, as required in Corollary \ref{structure_constants}).
\end{Remark}

\medskip

\section{A new proof of the P-R-V property for sums of minuscules} \label{PRV_section}

Before it was established independently by O. Mathieu \cite{Ma} and S. Kumar \cite{Ku}, the following was known as the P-R-V conjecture (see also \cite{Li} for a short proof based on Littelmann's path model).

\begin{theorem} [P-R-V property] \label{PRV_for_Rep}
If $\lambda = w_1 \mu_1 + \cdots w_r \mu_r$, then $V_\lambda$ appears with multiplicity at least 1 in the tensor product $V_{\mu_1} \otimes \cdots \otimes V_{\mu_r}$.
\end{theorem}
It is actually much easier to establish the Hecke-algebra analogue of the P-R-V property. 

\begin{prop} \label{PRV_for_Hecke}
Under the same assumptions as above, the function $f_\lambda$ appears in $f_{\mu_1} * \cdots * f_{\mu_r}$ with non-zero coefficient.
\end{prop}

\begin{proof} 
Recall that ${\rm Hecke}(\mu_\bullet, \lambda)$ holds if and only if the variety $m_{\mu_\bullet}^{-1}(e_\lambda) \cap {\mathcal Q}_{\mu_\bullet}$ is non-empty (see Lemma \ref{c(q)} and \cite{H}, $\S 4$).

But the equality $\lambda = w_1 \mu_1 + \cdots + w_r \mu_r$ yields a point $L_\bullet$ in the intersection 
$m_{\mu_\bullet}^{-1}(e_\lambda) \cap {\mathcal Q}_{\mu_\bullet}$, given by 
$$
L_i = t_{w_1\mu_1 + \cdots + w_i \mu_i} e_0,
$$
for $0 \leq i \leq r$.
\end{proof} 

Note that Corollary \ref{structure_constants} and Proposition \ref{PRV_for_Hecke} combine to give a new proof of Theorem \ref{PRV_for_Rep}, in the case where each $\mu_i$ is a sum of minuscule coweights (in particular for the group 
${\rm GL}_n$).

\medskip

\section{A saturation theorem for sums of minuscules}

The following saturation property for ${\rm Hecke}(\mu_\bullet,\lambda)$ is due to M. Kapovich, B. Leeb, and J. Millson \cite{KLM}.

\begin{theorem} [\cite{KLM}] \label{Hecke_saturation}  For any split group $G$ over $\mathbb F_q$, there exists a positive integer $k_G$ given explicitly in terms of the root data for $G$, with the following property: for any tuple of dominant coweights $(\mu_\bullet,\lambda)$ satisfying $\sum_i \mu_i - \lambda \in Q^\vee$, and every positive integer $N$, we have
$$
{\rm Hecke}(N\mu_\bullet,N\lambda) \Rightarrow {\rm Hecke}(k_G\mu_\bullet, k_G\lambda).
$$
\end{theorem}
We call $k_G$ the {\em Hecke algebra saturation factor for $\,G$}.
It turns out that $k_{{\rm GL}_n} = 1$, so this result shows that the structure constants for the Hecke algebra have the strongest possible saturation property in the case of ${\rm GL}_n$. 

\medskip

Corollary \ref{structure_constants} and Theorem \ref{Hecke_saturation} combine to give the following saturation theorem.

\begin{theorem} [Saturation for sums of minuscules -- weak form] \label{saturation_for_sums_of_minuscules}
Suppose that at least $r-1$ of the weights $\mu_i$ of $\widehat{G}$ are sums of dominant minuscule weights, and suppose the sum $\sum_i \mu_i$ belongs to the lattice spanned by the roots of $\widehat{G}$.  Let $N$ be any positive integer.  Then 
$$
(V_{N\mu_1} \otimes \cdots \otimes V_{N\mu_r})^{\widehat{G}} \neq 0 ~ \Rightarrow ~ (V_{k_G\mu_1} \otimes \cdots \otimes V_{k_G \mu_r})^{\widehat{G}} \neq 0.
$$
\end{theorem}

In the case of ${\rm GL}_n$ this was proved in \cite{KLM}, providing a new proof of the saturation property for ${\rm GL}_n$. 

\begin{proof}
Without loss of generality, we may assume $\mu_1, \dots, \mu_{r-1}$ are sums of minuscules.  Recall that for any highest weight representation $V_\mu$, the contragredient $(V_\mu)^*$ is also irreducible, so that we can define a dominant coweight $\mu^*$ by the equality $(V_\mu)^* = V_{\mu^*}$.  Let $\mu'_\bullet = (\mu_1, \dots, \mu_{r-1})$.  Now the theorem follows from Theorem \ref{Hecke_saturation} and Corollary \ref{structure_constants}, and the equivalences
\begin{align*}
{\rm Hecke}(\mu_\bullet, 0) &\Leftrightarrow {\rm Hecke}(\mu'_\bullet, \mu_r^*) \\
{\rm Rep}(\mu_\bullet, 0) &\Leftrightarrow {\rm Rep}(\mu'_\bullet, \mu_r^*).
\end{align*}
\end{proof}
In fact it seems that a stronger implication will hold.  Although it might not be necessary, here we will assume that {\em all} the weights $\mu_i$ are sums of minuscules, to be consistent with computer checks we ran with LiE.  We expect that the saturation factor $k_G$ can be omitted in the above statement.  

\begin{conj} [Saturation for sums of minuscules -- strong form] \label{conjecture}  Suppose each weight $\mu_i$ of $\widehat{G}$ is a sum of dominant minuscule weights, and suppose $\sum_i \mu_i$ belongs to the root lattice.  Then 
$$
(V_{N\mu_1} \otimes \cdots \otimes V_{N\mu_r})^{\widehat{G}} \neq 0 ~ \Rightarrow ~ (V_{\mu_1} \otimes \cdots \otimes V_{\mu_r})^{\widehat{G}} \neq 0.
$$
\end{conj}

\medskip

When $k_G$ is small (e.g. $k_{{\rm GSp}_{2n}} = 2$) the conjecture seems to be only a minor strengthening of Theorem \ref{saturation_for_sums_of_minuscules}.  However for some exceptional groups $k_G$ is quite large (e.g. $k_{E_7} = 12$) and there the conjecture indicates that a substantial strengthening of Theorem \ref{saturation_for_sums_of_minuscules} should remain valid.  In any case, the conjecture ``explains'' to a certain extent the phenomenon of saturation for ${\rm GL}_n$ by placing it in a broader context.

We present the following evidence for Conjecture \ref{conjecture}.  Taking Corollary \ref{structure_constants} into account, the following theorem results immediately from a slightly more comprehensive Hecke-algebra analogue, proved in a joint appendix with M. Kapovich and J. Millson (Theorem \ref{saturation_for_Hecke_constants}).

\begin{theorem} \label{evidence_for_conjecture}
Suppose that $G_{ad}$ a product of simple groups of type $A, B,C, D,$ or $E_7$.  
Suppose each $\mu_i$ is a sum of minuscules and that $\sum_i \mu_i \in Q(\widehat{G})$.  Then 
$$
{\rm Rep}(N\mu_\bullet,0) \Rightarrow {\rm Rep}(\mu_\bullet, 0),
$$
provided we assume either of the following conditions:
\begin{enumerate}
\item [(i)]  All simple factors of $G_{ad}$ are of type $A,B,C,$ or $E_7$;
\item [(ii)] All simple factors of $G_{ad}$ are of type $A,B,C,D$, or $E_7$, and for each factor of type $D_{2n}$ (resp. $D_{2n+1}$) the projection of $\mu_i$ onto that factor is a multiple of a single minuscule weight (resp. a multiple of the minuscule weight $\varpi_1$).
\end{enumerate} 
\end{theorem}

\begin{Remark} \label{KM_refinement}  M. Kapovich and J. Millson have recently announced in \cite{KM} that the implication 
$$
{\rm Rep}(N\mu_\bullet,0) \Rightarrow {\rm Rep}(k^2_G\mu_\bullet,0)
$$
holds for {\em every} split semi-simple group $G$ over $k(\!(t)\!)$ and for {\em all} weights $\mu_\bullet$ (assuming of course $\sum_i \mu_i \in Q(\widehat{G})$).  Conjecture \ref{conjecture} above is in a sense ``orthogonal'' to this statement: instead of fixing a group and then asking what saturation factor will work for that group, we are asking whether for certain special classes of weights $\mu_\bullet$ (e.g. sums of minuscules for groups that possess them) the saturation factor of 1 is guaranteed to work.  
\end{Remark}

\subsubsection{Relation with the conjecture of Knutson-Tao}

The following conjecture of Knutson-Tao proposes a sufficient condition on weights $\mu_1,\mu_2,\lambda$ of a general semi-simple group to ensure a saturation theorem will hold.

\begin{conj} [\cite{KT}] \label{KT_conj}
Let $\widehat{G}$ be a connected semi-simple complex group, and suppose $(\mu_1,\mu_2,\lambda)$ are weights of a maximal torus $\widehat{T}$ such that $\mu_1 + \mu_2 + \lambda$ annihilates all elements $s \in \widehat{T}$ whose centralizer in $\widehat{G}$ is a semi-simple group.  Then for any positive integer $N$,\begin{equation} \label{eq:saturation}
(V_{N\mu_1} \otimes V_{N\mu_2} \otimes V_{N\lambda})^{\widehat{G}} \neq 0 \Rightarrow (V_{\mu_1} \otimes V_{\mu_2} \otimes V_\lambda)^{\widehat{G}} \neq 0.
\end{equation}
\end{conj}

Fix a connected semi-simple complex group $\widehat{G}$.
It is natural to ask how Conjectures \ref{conjecture} and \ref{KT_conj} are related: if we assume $\mu_1, \mu_2$ and $\lambda$ are sums of minuscules, does the Knutson-Tao conjecture then imply Conjecture \ref{conjecture}?  The answer to this question is no, as the following example demonstrates.

\medskip

\noindent {\bf Example.}  Let $\widehat{G} = {\rm Spin}(12)$, the simply-connected group of type $D_6$.  Suppose $\mu_1, \mu_2, \lambda$ are three weights of $\widehat{T}$ whose sum belongs to the root lattice (so the sum annihilates the center $Z(\widehat{G})$).  Conjecture \ref{conjecture} asserts that (\ref{eq:saturation}) holds provided that 
\begin{equation} \label{eq:my_condition}
\mu_1, \mu_2, \lambda \in \mathbb Z_{\geq 0}[\varpi_1, \varpi_5, \varpi_6],
\end{equation}
where we have labeled characters using the conventions of \cite{Bou}.  Henceforth let us assume condition (\ref{eq:my_condition}).  Now, Conjecture \ref{KT_conj} asserts that (\ref{eq:saturation}) holds provided $\mu_1 + \mu_2 +\lambda$ also annihilates certain elements.  Consider the element $s:= \varpi^\vee_3(e^{2\pi i/2}) \in \widehat{T}$, an element of order 2.  It is easy to check that ${\rm Cent}_{\widehat{G}}(s)$ is a semi-simple group.  Furthermore, it is clear that $\mu_1 + \mu_2 + \lambda$ annihilates $s$ if and only if 
\begin{equation} \label{eq:KT_condition}
\langle \mu_1 + \mu_2 + \lambda, \varpi^\vee_3 \rangle \in 2\mathbb Z.
\end{equation}

But this last condition can easily fail: take for example $\mu_1 = \varpi_6 $, $\mu_2 = 0$, and $\lambda = \varpi_6$, so that $\mu_1 + \mu_2 + \lambda = 2\varpi_6 = e_1+e_2+e_3+e_4+e_5+e_6$ and thus
$$
\langle \mu_1 + \mu_2 + \lambda, \varpi_3^\vee \rangle = 3.
$$

In other words, if $\mu_1, \mu_2, \lambda$ are sums of minuscules for ${\rm Spin}(12)$, the Knutson-Tao conjecture predicts at most the implication
$$
(V_{N\mu_1} \otimes V_{N\mu_2} \otimes V_{N\lambda})^{\widehat{G}} \neq 0 \Rightarrow (V_{2\mu_1} \otimes V_{2\mu_2} \otimes V_{2\lambda})^{\widehat{G}} \neq 0,
$$
whereas Conjecture \ref{conjecture} predicts the sharper statement (\ref{eq:saturation}).  For the example $\mu_1 = \lambda = \varpi_6, \mu_2 = 0$ above, this sharper statement is indeed correct (use that $V_{\varpi_6}$ is a self-contragredient representation).

\medskip

\section{Equidimensionality of (locally closed) partial Springer varieties for ${\rm GL}_n$} \label{Springer_section}

In this section we will use Proposition \ref{further_equidim} to deduce similar equidimensionality results for ``locally closed'' Springer varieties associated to a partial Springer resolution of the nilpotent cone for ${\rm GL}_n$.  We will also characterize those which are non-empty and express the number of irreducible components in terms of structure constants.  Finally, we describe the relation of these questions with the Springer correspondence.  For the most part, our notation closely parallels that of \cite{BM}.

\medskip

\subsection{Definitions and the equidimensionality property}

Let $V$ denote a $k$-vector space of dimension $n$, and let $(d_1, \dots, d_r)$ denote an ordered $r$-tuple of nonnegative integers such that $d_1 + \cdots + d_r = n$.  The $r$-tuple $d_\bullet$ determines a standard parabolic subgroup $P \subset {\rm GL}(V) = {\rm GL}_n$.  We consider the variety of partial flags of type $P$:
$$
\mathcal P = \{ V_\bullet = (V = V_0 \supset V_1 \supset \cdots \supset V_r = 0) ~ | ~ {\rm dim}\Big(\dfrac{V_{i-1}}{V_i}\Big) = d_i, \,\, 1 \leq i \leq r \}.
$$

Consider the Levi decomposition $P = LN$, where $N$ is the unipotent radical of $P$, and $L \cong {\rm GL}_{d_1} \times \cdots \times {\rm GL}_{d_r}$.

For a nilpotent endomorphism $g \in {\rm End}(V)$, let ${\mathcal P}_g$ denote the closed subvariety of 
${\mathcal P}$ consisting of partial flags $V_\bullet$ such that $g$ stabilizes each $V_i$.  This is the Springer fiber (over $g$) of the partial Springer resolution 
$$
\xi: \widetilde{\mathcal N}^{\mathcal P} \rightarrow {\mathcal N}
$$
where ${\mathcal N} \subset {\rm End}(V)$ is the nilpotent cone, $\widetilde{\mathcal N}^{\mathcal P} = \{ (g, V_\bullet) \in {\mathcal N} \times {\mathcal P} ~ | ~ V_\bullet \in {\mathcal P}_g \}$, and the morphism $\xi$ forgets $V_\bullet$.   

The nilpotent cone ${\mathcal N}$ has a natural stratification indexed by the partitions of $n$.  These can be identified with dominant coweights $\lambda = (\lambda_1, \dots, \lambda_n)$ where $\lambda_1 \geq \cdots \geq \lambda_n \geq 0$ and $\lambda_1 + \cdots + \lambda_n = n$.  The integers $\lambda_i$ give the sizes of Jordan blocks in the normal form of an element in ${\mathcal N}$.  In a similar way, the partial Springer 
resolution $\widetilde{\mathcal N}^{\mathcal P}$ carries a natural stratification indexed by $r$-tuples 
$\mu'_\bullet = (\mu'_1, \dots, \mu'_r)$ where $\mu'_i$ is a partition of $d_i$ having length $n$ (see \cite{BM}, $\S 2.10$).  In other words, if we let 
$$
\mu_i = (d_i,0^{n-1})
$$ 
for $1 \leq i \leq r$, then $\widetilde{\mathcal N}^{\mathcal P}$ carries a natural stratification indexed by $r$-tuples $\mu'_\bullet = (\mu'_1, \dots, \mu'_r)$ where for each $i$, $\mu'_i$ is a dominant coweight for ${\rm GL}_n$ and $\mu'_i \preceq \mu_i$. The stratum indexed by $\mu'_\bullet$ consists of pairs $(g,V_\bullet)$ such that the Jordan form of the endomorphism on $V_{i-1}/V_i$ induced by $g$ has Jordan type $\mu'_i$, for $1 \leq i \leq r$. 

Let us denote by ${\mathcal N}_\lambda \subset {\mathcal N}$ the stratum indexed by $\lambda$.  Write $x = \mu'_\bullet$ for short and denote by ${\mathcal P}^{(x)}$ the stratum of $\widetilde{\mathcal N}^{\mathcal P}$ which is indexed by $x = \mu'_\bullet$.

The morphism $\xi: \widetilde{\mathcal N}^{\mathcal P} \rightarrow {\mathcal N}$ is a locally-trivial semi-small morphism of stratified spaces.  In fact, by the lemma below, it comes by restriction of the morphism $m_{\mu_\bullet} : \widetilde{\mathcal Q}_{\mu_\bullet} \rightarrow \overline{\mathcal Q}_{|\mu_\bullet|} = \overline{\mathcal Q}_{(n,0^{n-1})}$, where $\mu_i = (d_i, 0^{n-1})$, $1 \leq i \leq r$.

The following useful relation between the nilpotent cone and the affine Schubert variety $\overline{\mathcal Q}_{(n,0^{n-1})}$ was observed by Lusztig \cite{Lu1} and Ng\^{o} \cite{Ng}.  Here, the standard lattice $V \otimes_k \mathcal O \cong \mathcal O^n$ is the base-point in ${\mathcal Q}$, which we previously denoted by $e_0$.

\begin{lemma} \label{xi_is_restriction}

The morphism $g \mapsto (g + tI_n){\mathcal O}^n$ determines an open immersion $\iota: {\mathcal N} \hookrightarrow \overline{\mathcal Q}_{|\mu_\bullet|}$.  Furthermore, the restriction of $m_{\mu_\bullet}$ over ${\mathcal N}$ can be identified with $\xi$.  In other words, there is a Cartesian diagram
$$
\xymatrix{
\widetilde{\mathcal N}^{\mathcal P} \ar[r]^{\widetilde{\iota}} \ar[d]_{\xi} & \widetilde{\mathcal Q}_{\mu_\bullet} 
\ar[d]^{m_{\mu_\bullet}} \\ 
{\mathcal N} \ar[r]^{\iota} & \overline{\mathcal Q}_{|\mu_\bullet|}.}
$$

Moreover, for each $\lambda$ (resp. $\mu'_\bullet$), we have $\iota^{-1}({\mathcal Q}_\lambda) = {\mathcal N}_\lambda$ (resp. $\widetilde{\iota}^{-1}({\mathcal Q}_{\mu'_\bullet}) = {\mathcal P}^{(x)}$, where $x = \mu'_\bullet$).  
\end{lemma}

\begin{proof}
The fact that $g \mapsto (g + tI_n){\mathcal O}^n$ determines an open immersion is most easily justified by proving the analogous {\em global} statement.  We refer to the proof of \cite{Ng}, Lemme 2.2.2. for the proof, since this point is not crucial in our applications of this lemma.

The compatibility between $m_{\mu_\bullet}$ and $\xi$ can be found in that same paper of Ng\^{o} (he proves in loc. cit. Lemme 2.3.1 a corresponding global statement).  Since this compatibility is used below, we will sketch the proof.  It is a direct consquence of the following explicit description of the morphism $\widetilde{\iota}$.  

Suppose $(g,V_\bullet) \in \widetilde{\mathcal N}^{\mathcal P}$.  Since ${\rm deg}_t({\rm det}(g+tI_n)) = n$, the lattice $(g+tI_n)\mathcal O^n$ has $k$-codimension $n$ in $\mathcal O^n$.  So, we can identify the $k$-vector space $\mathcal O^n/(g + tI_n)\mathcal O^n$ with $V$, equivariantly for the action of $g \in 
{\rm End}(V)$ on both sides.  The $g$-stable partial flag $V_\bullet$ then determines a $g$-stable (hence also $t$-stable) partial flag in $\mathcal O^n/(g + tI_n)\mathcal O^n$.  Thus, we get a sequence of $\mathcal O$-lattices $\mathcal O^n = L_0 \supset L_1 \supset \cdots \supset L_r = (g + tI_n)\mathcal O^n$, such that ${\rm dim}_k(L_{i-1}/L_i) = d_i$, for all $1 \leq i \leq r$.  Hence $L_\bullet \in \widetilde{\mathcal Q}_{\mu_\bullet}$, and we have $\widetilde{\iota}(g,V_\bullet) = L_\bullet$.
\end{proof}

\medskip

The goal of this subsection is to prove, in Proposition \ref{equidim_for_Springer_varieties} below, an equidimensionality property of the {\em locally closed Springer varieties} $\, {\mathcal P}^{(x)}_y$.  We define these as follows.  Let $\lambda$ index the stratum ${\mathcal N}_\lambda$ of ${\mathcal N}$, and let $x = \mu'_\bullet$ index the stratum ${\mathcal P}^{(x)}$ of $\widetilde{\mathcal N}^{\mathcal P}$.  Let $y \in {\mathcal N}_\lambda$.  We define
$$
{\mathcal P}^{(x)}_y := \xi^{-1}(y) \cap {\mathcal P}^{(x)}.
$$ 
Put another way, 
\begin{equation} \label{eq:restriction_to_cone}
{\mathcal P}^{(x)}_y = {\widetilde{\iota}}^{-1}(m_{\mu_\bullet}^{-1}(y) \cap {\mathcal Q}_{\mu'_\bullet}).
\end{equation}
Further, let ${\mathcal P}^{x} = \overline{{\mathcal P}^{(x)}}$ and put ${\mathcal P}^x_y = \xi^{-1}(y) \cap {\mathcal P}^{x}$.  Thus,
\begin{equation}
{\mathcal P}^x_y = {\widetilde{\iota}}^{-1}(m_{\mu_\bullet}^{-1}(y) \cap \overline{\mathcal Q}_{\mu'_\bullet}).
\end{equation}  
This is essentially the notation used in \cite{BM}, $\S 3.2$.  Following loc. cit., we recall that 
\begin{itemize}
\item the {\em Steinberg variety} ${\mathcal P}_y := \xi^{-1}(y)$ is the disjoint union of its Springer parts ${\mathcal P}^{(x)}_y$, 
\item ${\mathcal P}_y = {\mathcal P}^{x}_y$ if $x$ is ``regular'' (i.e., $\mu'_i = \mu_i$, for all $1 \leq i \leq r$),
\item ${\mathcal P}^{(0)}_y = {\mathcal P}^0_y$ (where ``$x = 0$'' means $\mu'_i = \mu'_{i}(0) := (1^{d_i}, 0^{n - d_i})$ 
for all $1 \leq i \leq r$).
\end{itemize}

The varieties ${\mathcal P}^0_y$ are called {\em Spaltenstein} varieties in \cite{BM} and {\em Spaltenstein-Springer} varieties in \cite{H}. 

Now Proposition \ref{further_equidim} and (\ref{eq:restriction_to_cone}) immediately give us the following equidimensionality result for the locally-closed Springer varieties ${\mathcal P}^{(x)}_y$, where $x = \mu'_\bullet$ and $y \in {\mathcal N}_\lambda$.  It is quite possible that this result is already known to some experts, but it does not seem to appear in the literature.  In any case, the present proof via Proposition \ref{further_equidim}(2) is a very transparent one.

\begin{prop} \label{equidim_for_Springer_varieties}  
If ${\mathcal P}^{(x)}_y \neq \emptyset$, then every irreducible component of ${\mathcal P}^{(x)}_y$ has dimension $\langle \rho, |\mu'_\bullet| - \lambda \rangle$.

In particular, the Spaltenstein-Springer variety ${\mathcal P}^0_y$ (if non-empty), is equidimensional of dimension $\langle \rho, |\mu'_\bullet(0)| - \lambda \rangle$, where $\mu'_i(0) := (1^{d_i},0^{n-d_i})$ for each $i$.
\end{prop}

Note that the last statement was proved in \cite{Sp} [final Corollary], by completely different methods.  Spaltenstein also proved that the varieties ${\mathcal P}^0_y$ admit pavings by affine spaces, and this fact can now be seen as a special case of Corollary \ref{paved_by_affines}.

\begin{Remark} \label{nonequidim}
 Note that we have not proved the equidimensionality of the varieties ${\mathcal P}^x_y$, and indeed they are not always equidimensional.  In fact, it is known that there exist coweights $\mu_i = (d_i,0^{n-1})$ and $\lambda \prec |\mu_\bullet| = (n,0^{n-1})$ such that the partial Springer fiber ${\mathcal P}_y^x$ is {\em not} equidimensional for $y \in {\mathcal N}_\lambda$.  See \cite{St}, proof of Cor. 5.6, or \cite{Sh}, Thm. 4.15.  By virtue of Lemma \ref{xi_is_restriction}, the corresponding fibers $m_{\mu_\bullet}^{-1}(y)$ are also not equidimensional. 
\end{Remark}

\subsection{When are locally closed Springer varieties non-empty?}

Let $x = \mu'_\bullet$ and $y \in {\mathcal N}_\lambda$.  
It is clear that ${\mathcal P}^x_y \neq \emptyset$ if and only if $\lambda \preceq |\mu'_\bullet|$.  The non-emptiness of the locally closed varieties ${\mathcal P}^{(x)}_y$ is more subtle.

\begin{prop} \label{Springer_non-emptiness}
The locally-closed Springer variety ${\mathcal P}^{(x)}_y$ is non-empty if and only if ${\rm Rep}(\mu'_\bullet, \lambda)$ holds.  Furthermore, there are equalities 
\begin{align*}
\# \, \mbox{irreducible components of ${\mathcal P}^{(x)}_y$} &= \# \, \mbox{top-dimensional irreducible components of ${\mathcal P}^x_y$} \\ 
&= \# \, \mbox{top-dimensional irreducible components of $m_{\mu'_\bullet}^{-1}(y)$} \\
&= {\rm dim}(V^\lambda_{\mu'_\bullet}).
\end{align*}
\end{prop}

This is obvious from (\ref{eq:restriction_to_cone}) and our previous discussion.

\subsection{Relation with the Springer correspondence}

The question of when ${\mathcal P}^{(x)}_y \neq \emptyset$ can also be related to the Springer correspondence.  For this discussion we assume $k = \mathbb C$ and temporarily replace ${\rm GL}_n$ with any connected reductive group $G$.  The Springer correspondence is a cohomological realization of a one-to-one correspondence 
$$
\rho(y, \psi) \leftrightarrow (y, \psi)
$$
between irreducible representations $\rho$ of $W$ and the set of relevant pairs $(y,\psi)$, where $y$ is a stratum of ${\mathcal N}$ and $\psi$ is a representation of the fundamental group of that stratum, giving rise to a local system $L_\psi$.  See \cite{BM}, Theorem 2.2.  

Let $V_{(y,\psi)}$ denote the underlying vector space for the representation $\rho(y, \psi)$.  Then 
the Weyl group $W$ acts on the cohomology of the Steinberg variety
$$
{\rm H}^i({\mathcal B}_y, \mathbb Q), 
$$
and in fact if we let $d_y := {\rm dim}({\mathcal B}_y)$, we have the isomorphism of $W$-modules 
$$
{\rm H}^{2d_y}({\mathcal B}_y, \mathbb Q)^{\rho(y,1)} = V_{(y,1)},
$$
where the left-hand side denotes the isotypical component of type $\rho(y,1)$.  See \cite{BM}, $\S2.2$.

Now once again we assume $G = {\rm GL}_n$ (for the rest of this section).  
In this case, it is known that only the representations $\rho(y,1)$ arise, and they give a complete list of the irreducible representations of $W = S_n$.

In the sequel, the symbol $y$ will either denote a point $y \in {\mathcal N}_\lambda$, or the stratum $y = \lambda$ itself.  Similarly, sometimes $x$ will denote a point $x \in {\mathcal Q}_{\mu'_\bullet}$, and other times it will denote the stratum $x = \mu'_\bullet$ itself.  Hopefully context will make it clear what is meant in each 
case.  Note that $d_y = {\rm dim}({\mathcal B}_y) = \langle \rho, |\mu_\bullet| - \lambda \rangle$ in this case.

Let $W(L) = N_L(T)/T$ denote the Weyl group of the standard Levi subgroup $L$ of $P$ we already fixed.  Let ${\mathcal B}(L)$ (resp. ${\mathcal N}(L)$ ) denote the flag variety (resp. nilpotent cone) for $L$, and for $\ell \in {\mathcal N}(L)$, let ${\mathcal B}(L)_\ell$ denote the corresponding Steinberg variety.
As in \cite{BM}, $\S 2.10$, we can regard any index $x = \mu'_\bullet$ as corresponding to a unique nilpotent orbit $\ell$: the choice of $x$ and $\ell$ both amount to choosing an $r$-tuple $(\mu'_1, \dotsb, \mu'_r)$ where $\mu'_i$ is a partition of $d_i$ of length $n$.  Thus, we can also write ${\mathcal B}(L)_\ell = {\mathcal B}(L)_x$. 

The question of whether ${\mathcal P}^{(x)}_y \neq \emptyset$ is essentially equivalent to whether $\rho(x,1)$ appears in the restriction to $W(L)$ of the $W$-module ${\rm H}^{2d_y}({\mathcal B}_y,\mathbb Q)$.

\begin{prop} \label{restriction_to_Levi}
The locally-closed Springer variety ${\mathcal P}^{(x)}_y$ is non-empty if and only if the representation $\rho(x,1)$ of $W(L)$ appears with positive multiplicity in ${\rm H}^{2d_y}({\mathcal B}_y, \mathbb Q)|_{W(L)}$.  
Furthermore, the multiplicity is given by the formula
\begin{align*}
{\rm dim}_{\mathbb C}\, [{\rm Hom}_{W(L)}(\rho(x,1),{\rm H}^{2d_y}({\mathcal B}_y,\mathbb Q)|_{W(L)})] &= \# \, \mbox{top-dimensional irreducible components of ${\mathcal P}^x_y$} \\ 
&= {\rm dim}(V^\lambda_{\mu'_\bullet}).
\end{align*}
\end{prop}

\begin{proof}
Since $V_{(x,1)} = {\rm H}^{2d_x}({\mathcal B}(L)_x, \mathbb Q)^{\rho(x,1)}$, the first statement will follow from the proof of \cite{BM}, Theorem 3.3, which shows in effect that there is an isomorphism of $W(L)$-modules
$$
{\rm H}^{2d_y -2d_x}({\mathcal P}^x_y, {\rm IC}({\mathcal P}^x)) ~ \otimes ~ {\rm H}^{2d_x}({\mathcal B}(L)_x, \mathbb Q)^{\rho(x,1)} = [{\rm H}^{2d_y}({\mathcal B}_y, \mathbb Q)|_{W(L)}]^{\rho(x,1)}.
$$
Here ${\mathcal P}^x := \overline{{\mathcal P}^{(x)}}$ and ${\rm IC}({\mathcal P}^x)$ denotes the intersection complex of ${\mathcal P}^x$, following the conventions of loc. cit. (it is a complex supported in cohomological degrees $[0,{\rm dim}({\mathcal P}^x))$).

Now we note that, provided ${\mathcal P}^{(x)}_y \neq \emptyset$, 
$${\rm dim}({\mathcal P}^x_y) = \langle \rho, |\mu'_\bullet| - \lambda \rangle = \langle \rho, |\mu_\bullet| - \lambda \rangle - \langle \rho, |\mu_\bullet| - |\mu'_\bullet| \rangle = d_y - d_x.
$$ 
Here, we have used the isomorphism
$$
{\mathcal B}(L)_x \cong \eta^{-1}(x)
$$
of \cite{BM}, Lemma 2.10 (b) to justify the equality $d_x = \langle \rho, |\mu_\bullet| - |\mu'_\bullet| \rangle$.
Finally, it is well-known that since $\xi$ is semi-small, the dimension of ${\rm H}^{2d_y-2d_x}({\mathcal P}^x_y, {\rm IC}({\mathcal P}^x))$ is the number of irreducible components of ${\mathcal P}^x_y$ having dimension $d_y-d_x$ (see e.g. \cite{H}, Lemma 3.2).  By Proposition \ref{Springer_non-emptiness}, we are done.

\end{proof}

This gives a refinement and new proof of \cite{BM} Corollary 3.5, in the case of ${\rm GL}_n$.


\bigskip

\section{Appendix: constructing special $r$-gons in Bruhat-Tits buildings \\ 
by reduction to rank 1}

\begin{center} 
by Thomas J. Haines \and Michael Kapovich \and John J. Millson
\end{center}


\medskip

\subsection{Constructing $r$-gons with allowed side-lengths}

Let $G$ denote a connected reductive group over an algebraically closed field $k$.  Let $\mathcal O = k[\![t]\!]$ and $L = k (\!(t)\!)$. (We use the symbol $L$ in place of $F$ to emphasize that $k$ can be any algebraically closed field, and not just $\overline{\mathbb F}_p$ as in the main body of the paper.)  Further, let $K = G(\mathcal O)$ and define the affine Grassmannian ${\mathcal Q} = G(L)/K$, viewed as an ind-scheme over $k$.  We fix once and for all a maximal torus $T \subset G$ and a Borel subgroup $B = TU$ containing $T$.

For a cocharacter $\mu$ of $T$, we shall denote by $\overline{\mu}$ the cocharacter of the adjoint group $G_{ad}$ which results by composing $\mu$ with the homomorphism $G \rightarrow G_{ad}$.  Recall that $G_{ad}$ is a product of simple adjoint groups $H$, and we will denote by $\overline{\mu}_H$ the composition of $\overline{\mu}$ with the projection $G_{ad} \rightarrow H$.   We have $\overline{\mu}_H \in X_*(T_H)$, where $T_H$ is the image of the torus $T$ under the homomorphism $G \rightarrow H$. 

Throughout this appendix, {\em dominant coweight} means $B$-{\em dominant cocharacter}.  Similar terminology will apply to the quotients $H$ (we use the Borel $B_H$ which is the image of $B$).

Recall that each factor $H$ corresponds to an irreducible finite root system whose Weyl group possesses a unique longest element $w_{H,0}$.  For any coweight $\nu$ of $T_H$, we set $\nu^* = -w_{H,0}\nu$.  We call such a coweight $\nu$ {\em self-dual} if $\nu^* = \nu$.  


Let $Q^\vee(H)$ (resp. $P^\vee(H) = X_*(T_H)$) denote the coroot (resp. coweight) lattice of the adjoint group $H$.

Our first result is the following generalization of Proposition 7.7 of \cite{KM}.  To state it, we need to single out a special class of fundamental coweights.

\begin{definition} Let $\varpi^\vee_i$ denote a a fundamental coweight of an adjoint group $H$.  We call $\varpi^\vee_i$ {\em allowed} if it satisfies the following properties: 
\begin{enumerate}
\item $\varpi^\vee_i$ is self-dual;
\item $n \varpi^\vee_i \in Q^\vee(H) \Leftrightarrow n \in 2 \, \mathbb Z$.
\end{enumerate}
\end{definition}

\begin{prop} \label{KM_generalized}
Suppose that for each factor $H$ of $G_{ad}$ we are given an allowed fundamental 
coweight 
$\lambda_H \in X_*(T_H)$. 

Suppose that for each $i = 1, 2, \dots, r$, the image $\overline{\mu_i}$ of the dominant coweight  $\mu_i \in X_*(T)$ is 
a sum of the form
$$
\overline{\mu_i} = \sum_H a^H_i \lambda_H
$$
for nonnegative integers $a^H_i$.  Suppose that $\sum_i \mu_i \in Q^\vee(G)$ and that the coweights ${\mu}_i$ satisfy the following {\em weak generalized triangle inequalities}
\begin{equation} 
{\mu}^*_i \preceq {\mu}_1 + \cdots + \widehat{{\mu}_i} + \cdots + {\mu}_r. \label{eq:weak_triangle_inequalities}
\end{equation}

Then the variety ${\mathcal Q}_{\mu_\bullet} \cap m_{\mu_\bullet}^{-1}(e_0)$ is non-empty.
\end{prop}

In the terminology of \cite{KLM}, \cite{KM}, the building of $G(L)$ has a closed $r$-gon with side-lengths $\mu_1, \dots, \mu_r$, whose vertices are special vertices.  We call these {\em special $r$-gons}.

\begin{proof}
For each factor $H$, let $\alpha_H$ denote the simple $B$-positive root corresponding to the fundamental coweight $\lambda_H$.  We consider the Levi subgroup $M \subset G$ that is generated by $T$ along with the root groups for all the roots $\pm \alpha_H$:
$$
M := \langle T, \,\, U_{\alpha_H}, \, U_{-\alpha_H} \, \rangle.
$$
The coweights $\mu_i$ resp. $\lambda_H$ determine coweights for $M$ resp. $M_{ad}$; we write $\overline{\overline{\mu_i}}$ resp. $\overline{\overline{\lambda_H}}$ for their images in the adjoint group $M_{ad}$.  Note that 
$$
M_{ad} = \prod_H {\rm PGL}_2,
$$
and that in the factor indexed by $H$, we can identify $\overline{\overline{\alpha_H}} = e_1 - e_2$ and $\overline{\overline{\lambda_H}} = (1,0)$. 

Now our assumptions imply that for each factor $H$, 
\begin{itemize}
\item $a^H_1 + \cdots + a^H_r$ is even, and 
\item $a^H_i \leq a^H_1 + \cdots + \widehat{a^H_i} + \cdots + a^H_r$, for each $i$.
\end{itemize}

As we shall see in the next lemma, these properties imply that there is a special $r$-gon in the building for ${\rm PGL}_2$ with side-lengths $a^H_1, \dots, a^H_r$.  Note that since $k$ is infinite, we will be working with a tree having infinite valence at each vertex, but this causes no problems.

\begin{lemma} \label{triangles_in_trees}
Suppose $u_1, \dots, u_r$ are nonnegative integers satisfying the generalized triangle inequalities
\begin{align*}
u_i &\leq u_1 + \cdots + \widehat{u_i} + \cdots + u_r, \,\,\, \forall i, \\
\sum_i u_i &\equiv 0 \,\, (\mbox{mod $2$}).
\end{align*}
Then there exists a special $r$-gon in the tree ${\mathcal B}({\rm PGL}_2)$ having side lengths $u_1, \dots, u_r$. 
\end{lemma}

\begin{proof}
First we claim that there exist integers $l$ and $m$, with $1 \leq l < m \leq r$, such that if we set 

\begin{align*}
A &= u_1 + \cdots + u_l \\
B &= u_{l+1} + \cdots + u_m \\
C &= u_{m+1} + \cdots + u_r,
\end{align*}

then 
\begin{align*}
A &\leq B + C \\
B &\leq A + C \\
C &\leq A + B\\
A+B+C &\equiv 0 \,\,\, (\mbox{mod $2$}).
\end{align*}
Indeed, note that for all $i$, $u_i \leq \frac{1}{2}(\sum_i u_i)$. We may choose $l$ to be the largest such that 
$$
u_1 + \cdots + u_{l} \leq \frac{1}{2}(\sum_i u_i),
$$
and then set $m = l+1$ (note that necessarily $l \leq r-1$, if at least one  $u_i > 0$).  

Now given $A,B,C$ as above, we may construct a ``tripod'' in the building as follows.  Choose any vertex $v_0$, and construct a tripod, centered at $v_0$, with legs having lengths $l_1, l_2, l_3$, where
\begin{align*}
l_1 &= \dfrac{A + B - C}{2} \\
l_2 &= \dfrac{B + C - A}{2} \\
l_3 &= \dfrac{A + C - B}{2}.
\end{align*}

Let $v_1,v_2,v_3$ denote the extreme points of the tripod, and consider the oriented paths $[v_0,v_i]$, where we have labeled vertices in such a way the length of $[v_0,v_i]$ is $l_i$. 

This yields a special 3-gon: the three ``sides'' are 
\begin{align*} 
[v_3,v_0] &\cup [v_0,v_1], \,\,\,\, (\mbox{length $=A$}) \\
[v_1,v_0] &\cup [v_0,v_2], \,\,\,\, (\mbox{length $=B$}) \\
[v_2,v_0] &\cup [v_0,v_3], \,\,\,\, (\mbox{length $=C$}).
\end{align*}

This (oriented) triangle begins and ends at the special vertex $v_3$.  The sides are themselves partitioned into smaller intervals of lengths $u_1, u_2, \dots, u_l$, etc.  Thus we have a special $r$-gon with the desired side-lengths.
\end{proof}

The building for $M_{ad}$ is simply the direct product of the buildings for the various ${\rm PGL}_2$ factors, hence we have a special $r$-gon in the building for $M_{ad}$ with side lengths $\overline{\overline{\mu}}_1, \dots, \overline{\overline{\mu}}_r$. Equivalently, we have
$$
1_{M_{ad}} \in M_{ad}(\mathcal O) \overline{\overline{\mu}}_1 M_{ad}(\mathcal O) \cdots M_{ad}(\mathcal O) \overline{\overline{\mu}}_r M_{ad}(\mathcal O).
$$
Now we want to claim that this implies that 
$$
1_M \in M(\mathcal O) \mu_1 M(\mathcal O) \cdots M(\mathcal O) \mu_r M(\mathcal O).
$$
But this follows from the Lemma \ref{reduction} below.  Since $M(\mathcal O) \subset K$, this immediately implies that 
$$
1_G \in K \mu_1 K \cdots K \mu_r K,
$$
and thus ${\mathcal Q}_{\mu_\bullet} \cap m_{\mu_\bullet}^{-1}(e_0) \neq \emptyset$, as desired.
\end{proof}

\subsection{Enumerating allowed coweights in $H$}  
Proposition \ref{KM_generalized} is most interesting for groups which are not of type $A$.  
For each type of adjoint simple factor $H$ not of type $A$, we enumerate the allowed and the minuscule fundamental coweights.  
We follow the indexing conventions of \cite{Bou} (note that our coweights are weights for the dual root system). 

\nopagebreak

\vskip.3cm

{\small\sf
\begin{tabular}{|c|c|c|}
\hline
Type of $H$ & Allowed fundamental coweights & Minuscule coweights \\
\hline \hline

$B_n$ & $\varpi^\vee_i, \,\, \mbox{$i$ odd}$ & $\varpi^\vee_1$ \\ \hline
$C_n$ & $\varpi^\vee_n$ & $\varpi^\vee_n$ \\ \hline
$D_{2n}$ & $\varpi^\vee_{2i - 1} \,\,\, (1 \leq i \leq n-1), \,\,\, \varpi^\vee_{2n-1}, \varpi^\vee_{2n}$ & 
$\varpi^\vee_1, \varpi^\vee_{2n-1}, \varpi^\vee_{2n}$  \\ \hline
$D_{2n+1}$ & $\varpi^\vee_{2i - 1}, \,\,\, 1 \leq i \leq n$ & $\varpi^\vee_1, \varpi^\vee_{2n}, 
\varpi^\vee_{2n+1}$ \\ \hline
$E_6$ & - & $\varpi^\vee_1, \varpi^\vee_6$ \\ \hline
$E_7$ & $\varpi^\vee_2, \varpi^\vee_5, \varpi^\vee_7$ & $\varpi^\vee_7$ \\ \hline
$E_8$ & - & - \\ \hline
$F_4$ & - & - \\ \hline
$G_2$ & - & - \\\hline
\end{tabular}
}
\vskip.3cm

Note that for each case where $H$ possesses a unique minuscule coweight ($B_n$, $C_n$, $E_7$), that minuscule coweight is allowed.  For type $D_{2n}$, all three minuscule coweights are allowed, but for $D_{2n+1}$, only $\varpi^\vee_1$ is allowed. 

\bigskip

\subsection{Reduction to adjoint groups}

\begin{lemma} \label{reduction}
For any connected reductive algebraic group $G$ over $k$ and dominant coweights $\mu_i$ such that $\sum_i \mu_i \in Q^\vee(G)$, we have the following two statements.
\begin{enumerate}
\item [(1)]  The canonical homomorphism $\pi: G(\mathcal O) \rightarrow G_{ad}(\mathcal O)$ is surjective.
\item [(2)]  Let $Z$ denote the center of $G$.  If $z \in Z(L)$ satisfies
$$
z1_G  \in G(\mathcal O) \mu_1 G(\mathcal O) \cdots G(\mathcal O) \mu_r G(\mathcal O)
$$
then $z \in Z(\mathcal O) \subset G(\mathcal O)$, and thus the statement holds with $z$ omitted.
\end{enumerate} 
\end{lemma}

\begin{proof}
Property (1) is a standard fact resulting from Hensel's lemma (see e.g. \cite{PR}, Lemma 6.5).  Let us recall briefly the proof.  For $a \in G_{ad}(\mathcal O)$, the preimage $\pi^{-1}(a)$ in $G(\mathcal O)$ is the set of $\mathcal O$-points of a smooth $\mathcal O$-scheme (a torsor for the smooth $\mathcal O$-group scheme $Z_{\mathcal O}$).  Clearly the reduction modulo $t$ of $\pi^{-1}(a)$ has a $k$-point (the residue field $\mathcal O/(t) = k$ being assumed algebraically closed).  Now by Hensel's lemma, $\pi^{-1}(a)$ also has an $\mathcal O$-point, proving property (1).

For Property (2) let us consider first the case of ${\rm GL}_n$.  The hypothesis implies that $zI_n$ belongs to the kernel of the homomorphism
$$
{\rm val} \circ {\rm det} : GL_n(L) \rightarrow \mathbb Z,
$$
since both $K$ and $Q^\vee({\rm GL}_n) \hookrightarrow T(L)$ belong to that kernel.  But then it is clear that $z \in \mathcal O^\times$.

In general, the same argument works if we replace ${\rm val} \circ {\rm det}$ with the Kottwitz homomorphism
$$
\omega_G : G(L) \rightarrow X^*(Z(\widehat{G}))_I,
$$
where $I$ denotes the inertia group ${\rm Gal}(L^{sep}/L)$ \footnote{Since $G$ is split over $L$, we may omit the coinvariants under $I$ here.};  see \cite{Ko}, $\S 7$ for the construction and properties of this map (we will use the functoriality of $G \mapsto \omega_G$ below).  Indeed, since $G(\mathcal O)$ and $Q^\vee(G) \hookrightarrow T(L)$ belong to the kernel of $\omega_G$, so does $z$.  Therefore, we will be done once we justify the equality  
$$Z(L) \cap {\rm ker}(\omega_G) = Z(\mathcal O).
$$   
Suppose $z \in Z(L) \cap {\rm ker}(\omega_G)$.  Let $T$ denote a (split) maximal $k$-torus of 
$G$, and consider the composition 
$$
\xymatrix{
Z \ar[r] & T \ar[r]^{c \,\,\,\,\,\,\,\,\,\,\,\,\,\,} & D := G/G_{der}.
}
$$  
By the functoriality of $\omega_G$, $z \in {\rm ker}(\omega_G)$ implies that $c(z) \in {\rm ker}(\omega_D)$.  Since $D_L$ is a split torus over $L$, the latter kernel is $D(\mathcal O)$.  Now since 
$$
Z(\mathcal O) \rightarrow D(\mathcal O)
$$
is surjective (by the same proof as in part (1)), there exists $z_0 \in Z(\mathcal O)$ such that 
$$
z_0^{-1}z \in {\rm ker}[G(L) \rightarrow D(L)] = G_{der}(L).
$$
But $Z(L) \cap G_{der}(L) = Z(k) \cap G_{der}(k)$ (since the latter is a finite group), which obviously belongs to $Z(\mathcal O)$.  This implies that $z \in Z(\mathcal O)$, as claimed.
\end{proof}

As a corollary of the proof, we have
\begin{cor} \label{reduction_to_ad}  For any tuple $(\mu_\bullet, \lambda)$ such that $\sum_i \mu_i - \lambda \in Q^\vee(G)$, 
$${\rm Hecke}^G(\mu_\bullet, \lambda) \Leftrightarrow {\rm Hecke}^{G_{ad}}(\overline{\mu}_\bullet, \overline{\lambda}).$$
\end{cor}

\medskip

The dual of the homomorphism 
$$
T \rightarrow T_{ad}
$$
is the composition
$$
\widehat{T}_{sc} \twoheadrightarrow \widehat{T}_{der} \hookrightarrow \widehat{T},
$$ 
where $\widehat{T}_{der} := \widehat{T} \cap \widehat{G}_{der}$ and $\widehat{T}_{sc}$ is the preimage of $\widehat{T}_{der}$ under the isogeny $\widehat{G}_{sc} \rightarrow \widehat{G}_{der}$.  Viewing a coweight $\lambda \in X_*(T)$ as a weight for the dual torus $\widehat{T}$, we let $\overline{\lambda}$ denote its image under the map
$$
X^*(\widehat{T}) \rightarrow X^*(\widehat{T}_{sc}).
$$

With this notation, Corollary \ref{reduction_to_ad} has the following analogue.  

\begin{lemma} \label{reduction_to_sc}   For any tuple $(\mu_\bullet, \lambda)$ of weights such that $\sum_i \mu_i - \lambda \in Q(\widehat{G})$, 
$${\rm Rep}^{\widehat{G}}(\mu_\bullet, \lambda) \Leftrightarrow {\rm Rep}^{\widehat{G}_{sc}}(\overline{\mu}_\bullet, \overline{\lambda}).$$
\end{lemma} 

\begin{proof}
Use the fact that the restriction of $V_{\lambda} \in {\rm Rep}(\widehat{G})$  along $\widehat{G}_{sc} \rightarrow \widehat{G}$ is simply $V_{\overline{\lambda}} \in {\rm Rep}(\widehat{G}_{sc})$.
\end{proof}

\bigskip

\subsection{A saturation theorem for Hecke algebra structure constants}

Assume now that $\lambda = 0$ and $\sum_i \mu_i \in Q^\vee(G)$.

\begin{theorem} \label{saturation_for_Hecke_constants}
Suppose that $G_{ad}$ a product of simple groups of type $A, B,C, D,$ or $E_7$.  Suppose that the projection of each $\mu_i$ onto a simple adjoint factor of $G_{ad}$ having type $B,C,D$ or $E_7$ is a multiple of a single allowed coweight.   Then 
$$
{\rm Hecke}(N\mu_\bullet,0) \Rightarrow {\rm Hecke}(\mu_\bullet, 0).
$$

In particular, this conclusion holds if each $\mu_i$ is a sum of minuscule coweights, provided we assume either of the following conditions:
\begin{enumerate}
\item [(i)]  All simple factors of $G_{ad}$ are of type $A,B,C,$ or $E_7$;
\item [(ii)] All simple factors of $G_{ad}$ are of type $A,B,C,D$, or $E_7$, and for each factor of type $D_{2n}$ (resp. $D_{2n+1}$) the projection of $\mu_i$ onto that factor is a multiple of a single minuscule coweight (resp. a multiple of the minuscule coweight $\varpi^\vee_1$).
\end{enumerate} 
\end{theorem}

\begin{proof}
By Corollary \ref{reduction_to_ad} we can assume $G$ is adjoint, and then prove the saturation property 
one factor at a time.  For factors of type $A$, the desired saturation property follows from \cite{KLM}, Theorem 1.8.  For factors of type $B,C,D$ or $E_7$, observe that the assumption ${\rm Hecke}(N\mu_\bullet,0)$ implies that the weak generalized triangle inequalities (\ref{eq:weak_triangle_inequalities}) hold, and then use Proposition \ref{KM_generalized}.
\end{proof}
 
When each $\mu_i$ is a sum of minuscules, it is very probable that the implication holds with no assumption on $G_{ad}$, in other words, factors $D$ and $E_6$ should be allowed in (i) (see Conjecture \ref{conjecture}).
There is ample computer evidence corroborating this.  However, the method of reduction to rank 1 used above breaks down for type $E_6$ and yields only limited information for type $D$, and thus a new idea seems to be required.

\section{Erratum}

\subsection{Correction to proof of Lemma 9.4(2)}

The proof of Lemma \ref{reduction} is flawed when ${\rm char}(k) >0 $ because $Z$ need not be smooth over ${\rm Spec}(k)$ (e.g. $G = {\rm SL}_n$ for  a field $k$ of characteristic dividing $n$).  In particular (1) is false in general, e.g. for $G = {\rm SL}_2$ when $k = \bar{\mathbb F}_2$. Nevertheless, part (2) is true as stated.

\subsubsection{Connected center case}
The proof of (1) works fine when $Z$ is connected (it is then a torus over $k$, hence smooth; hence $Z_{\mathcal O}$ is a smooth scheme over $\mathcal O$).  

The proof of (2) concerning the step $Z(L) \cap {\rm ker}(\omega_G) \subseteq Z(\mathcal O)$ needs to be altered, because even if $Z$ is connected, the map $Z(\mathcal O) \rightarrow D(\mathcal O)$ need not be surjective, since the kernel of $Z \rightarrow D$ is finite and might not be smooth.  However, we can argue as follows.  Since $Z$ is connected, there is a commutative diagram with exact rows 
$$
\xymatrix{
1 \ar[r] & Z(L) \ar[r] \ar[d]_{\omega_Z} & G(L) \ar[r] \ar[d]_{\omega_G} & G_{\rm ad}(L) \ar[r] \ar[d]_{\omega_{G_{\rm ad}}} & 1 \\
0 \ar[r] & \pi_1(Z) \ar[r] & \pi_1(G) \ar[r] & \pi_1(G_{\rm ad}) \ar[r] & 0.}
$$
Therefore $Z(L) \cap {\rm ker}(\omega_G) \subseteq Z(L) \cap {\rm ker}(\omega_Z) = Z(\mathcal O)$, as desired.

\subsubsection{Reduction to connected center case}
We need to show we can reduce (2) to the connected center case.

Given $G$ with center $Z$ and maximal torus $T$, we can embed $G \hookrightarrow \tilde{G}$ as a normal subgroup, where $\tilde{G}$ is a connected reductive group over $k$ whose center $\tilde{Z}$ contains $Z$ and is connected. Moreover $G_{\rm der} = \tilde{G}_{\rm der}$ (so the root systems for $G$ and $\tilde{G}$ coincide). Therefore $G_{\rm ad} = \tilde{G}_{\rm ad}$, and hence $\tilde{Z} \cap G = Z$. Explicitly, write $G = G_{\rm der} \times^{Z_{\rm der}} Z$ where the quotient group is formed using the antidiagonal action of $Z_{\rm der}$, and set $\tilde{G} = G_{\rm der} \times^{Z_{\rm der}} T$.  Then $\tilde{G}$ has center $\tilde{Z} = Z_{\rm der} \times^{Z_{\rm der}} T$ and maximal torus $\tilde{T} = T_{\rm der} \times^{Z_{\rm der}} T$.  By construction $\tilde{T} \cap G = T$.

Now suppose in proving (2) for $G$ that for $z \in Z(L) \subset \tilde{Z}(L)$
$$
z 1_G \in G(\mathcal O) \mu_1 \G(\mathcal O) \cdots G(\mathcal O)\mu_r  G(\mathcal O).
$$
Using $G(\mathcal O) \subset \tilde{G}(\mathcal O)$ and the result for $\tilde{G}$, we obtain $z \in G(L) \cap  \tilde{Z}(\mathcal O) = Z(\mathcal O)$.  Hence (2) is proved.

\subsection{Correction to proof of Corollary 9.5}

We give a correct proof of the implication ${\rm Hecke}^G(\mu_\bullet, \lambda) \Leftarrow {\rm Hecke}^{G_{ad}}(\overline{\mu}_\bullet, \overline{\lambda}).$ 

\begin{lemma} \label{O-surj}
The multiplication map $m: G \times \tilde{T} \rightarrow \tilde{G}$ is a smooth and surjective morphism of affine $k$-schemes, hence by Hensel's lemma induces a surjective map $G(\mathcal O) \times \tilde{T}(\mathcal O) \rightarrow \tilde{G}(\mathcal O)$.
\end{lemma}

\begin{proof}
Clearly $m$ is surjective, and since $\tilde{G}$ is an integral $k$-scheme, it is generically flat. The transitive action of $G \times \tilde{T}$ shows it is flat everywhere. Since every geometric fiber is isomorphic to $T$, $m$ is smooth.
\end{proof}

Suppose
$$
t^{\bar{\lambda}} \in G_{\rm ad}({\mathcal O})t^{\bar{\mu}_1}G_{\rm ad}(\mathcal O) \cdots G_{\rm ad}(\mathcal O) t^{\bar{\mu}_r} G_{\rm ad}(\mathcal O).
$$
By the connected center case of Lemma \ref{reduction}(1), there is some $\tilde{z} \in \tilde{Z}(L)$ such that
$$
t^{\lambda} \tilde{z} \in \tilde{G}({\mathcal O})t^{{\mu}_1}\tilde{G}(\mathcal O) \cdots \tilde{G}(\mathcal O) t^{{\mu}_r} \tilde{G}(\mathcal O).
$$
By Lemma \ref{reduction}(2), $\tilde{z} \in \tilde{Z}(\mathcal O)$, so we can omit it. Then using $\tilde{G}(\mathcal O) = G(\mathcal O)\tilde{T}(\mathcal O)$ and $G(\mathcal O) \triangleleft \tilde{G}(\mathcal O)$, we may write
$$
t^{\lambda} \tilde{t} \in G({\mathcal O})t^{{\mu}_1}G(\mathcal O) \cdots G(\mathcal O) t^{{\mu}_r} G(\mathcal O)
$$
for some $\tilde{t} \in \tilde{T}(\mathcal O)$. But then $\tilde{t} \in \tilde{T}(\mathcal O) \cap G(L) = T(\mathcal O)$, and so we can omit $\tilde{t}$ as well.
\qed

\bigskip
\bigskip

\obeylines
Mathematics Department
University of Maryland
College Park, MD 20742-4015
tjh@math.umd.edu

\end{document}